\documentclass[12pt]{article}
\usepackage{amsmath,amsfonts,amssymb,amsthm}
\begin{document}
\def \Z{\mathbb Z}
\def \C{\mathbb C}
\def \A{{\mathcal{A}}}
\def \D{{\mathcal{D}}}
\def \E{{\mathcal{E}}}
\def \H{{\mathcal{H}}}
\def \S{{\mathcal{S}}}
\def \R{\mathbb R}
\def \Q{\mathbb Q}
\def \N{\mathbb N}
\def \P{\mathbb P}
\def \bZ{{\bf Z}}
\def \wt{{\rm wt}\;}
\def \gr{{\rm gr}}
\def \g{\mathfrak{g}}
\def \gl{\mathfrak{gl}}
\def \sl{\mathfrak{sl}}
\def \h{\mathfrak h}
\def \mod{{\rm mod}\;}

\def \Aut{{\rm Aut}}
\def \GL{{\rm GL}}
\def \PAut{{\rm PAut}}
\def \Der{{\rm Der}}
\def \PDer{{\rm PDer}}
\def \Sing{{\rm Sing}}
\def \span{{\rm span}}
\def \Res{{\rm Res}}
\def \End{{\rm End}}
\def \PEnd{{\rm PEnd}}
\def \Hom{{\rm Hom}}
\def \<{\langle} 
\def \>{\rangle}
\def \be{\begin{equation}\label}
\def \ee{\end{equation}}
\def \bex{\begin{exa}\label}
\def \eex{\end{exa}}
\def \bl{\begin{lem}\label}
\def \el{\end{lem}}
\def \bt{\begin{thm}\label}
\def \et{\end{thm}}
\def \bp{\begin{prop}\label}
\def \ep{\end{prop}}
\def \br{\begin{rem}\label}
\def \er{\end{rem}}
\def \bc{\begin{coro}\label}
\def \ec{\end{coro}}
\def \bd{\begin{de}\label}
\def \ed{\end{de}}
\def \bconj{\begin{conj}\label}
\def \econj{\end{conj}}

\newtheorem{thm}{Theorem}[section]
\newtheorem{prop}[thm]{Proposition}
\newtheorem{coro}[thm]{Corollary}
\newtheorem{conj}[thm]{Conjecture}
\newtheorem{exa}[thm]{Example}
\newtheorem{lem}[thm]{Lemma}
\newtheorem{rem}[thm]{Remark}
\newtheorem{de}[thm]{Definition}
\newtheorem{hy}[thm]{Hypothesis}
\makeatletter
\@addtoreset{equation}{section}
\def\theequation{\thesection.\arabic{equation}}
\makeatother
\makeatletter

\begin{center}
{\Large \bf Constructing quantum vertex algebras}
\end{center}
\begin{center}{Haisheng Li\footnote{Partially supported by an NSA grant}\\
Department of Mathematical Sciences, Rutgers University, Camden, NJ 08102\\
and\\
Department of Mathematics, Harbin Normal University, Harbin, China}
\end{center}

\begin{abstract}
This is a sequel to \cite{li-qva}. In this paper, we focus on
the construction of quantum vertex algebras over $\C$,
whose notion was formulated in \cite{li-qva} with Etingof and Kazhdan's
notion of quantum vertex operator algebra (over $\C[[h]]$) as one of
the main motivations.  As one of the main steps in constructing
quantum vertex algebras, we prove that every countable-dimensional
nonlocal (namely noncommutative) vertex algebra over $\C$, which
either is irreducible or has a basis of PBW type, is nondegenerate in
the sense of Etingof and Kazhdan.  Using this result, we establish the
nondegeneracy of better known vertex operator algebras and some
newly constructed nonlocal vertex algebras.  We construct a family of quantum
vertex algebras closely related to Zamolodchikov-Faddeev algebras.
\end{abstract}

\section{Introduction}
In a previous paper \cite{li-qva}, to seek for appropriate vertex
algebra-like objects for quantum affine algebras we studied fairly general
formal vertex operators on an arbitrary vector space and we found that
general formal vertex operators such as the generating functions of
any quantum affine algebra in Drinfed's realization (see \cite{dr}), 
acting on a highest weight module (see \cite{fj}), 
did give rise to vertex algebra-like structures
in certain natural way.  
More specifically, given any vector space $W$ (over $\C$), we
studied what we called quasi compatible subsets of $\Hom (W,W((x)))$
and we proved that every maximal quasi compatible subspace has a
natural nonlocal (namely noncommutative) vertex algebra structure and
that any quasi compatible subset generates a nonlocal vertex algebra.
(Nonlocal vertex algebras were also called field
algebras in \cite{bk} and $G_{1}$-vertex algebras in \cite{li-g1}.)
We found that for a quasi compatible set
of formal vertex operators of a certain special type,
the generated nonlocal vertex algebra very much resembles
a quantum vertex operator algebra in the sense of Etingof and Kazhdan
(see \cite{ek}).
This naturally led us to certain notions of 
weak quantum vertex algebra and quantum vertex algebra over $\C$.

Nonlocal vertex algebras are both vertex analogues and generalizations of
noncommutative associative algebras and they are too general to work with.
Naturally, one seeks for a good family of nonlocal vertex algebras 
which have more interesting structures.
Then weak quantum vertex algebras are singled out for this consideration.
Weak quantum vertex algebras are nonlocal vertex algebras (which are vertex analogues of 
noncommutative associative algebras), satisfying a certain generalized locality (where
locality is a vertex analogue of classical commutativity).
(Quantum vertex operator algebras in the sense of Etingof and Kazhdan
are $h$-adic nonlocal vertex algebras satisfying 
a certain $h$-adic generalized locality.)
Precisely, a weak quantum vertex algebra 
can be defined by using all the axioms
in defining a vertex algebra $V$ except the Jacobi identity axiom
which is replaced by the following axiom:
For any $u,v\in V$, there exists
$\sum_{i=1}^{r}v^{(i)}\otimes u^{(i)}\otimes f_{i}(x)\in V\otimes V\otimes \C((x))$ such that
\begin{eqnarray}\label{e1.1}
& &x_{0}^{-1}\delta\left(\frac{x_{1}-x_{2}}{x_{0}}\right)Y(u,x_{1})Y(v,x_{2})
-x_{0}^{-1}\delta\left(\frac{x_{2}-x_{1}}{-x_{0}}\right)
\sum_{i=1}^{r}f_{i}(-x_{0})Y(v^{(i)},x_{2})Y(u^{(i)},x_{1})\nonumber\\
& &=x_{2}^{-1}\delta\left(\frac{x_{1}-x_{0}}{x_{2}}\right)
Y(Y(u,x_{0})v,x_{2}).
\end{eqnarray}
Clearly, this notion generalizes
the notions of vertex algebra, vertex superalgebra and vertex color algebra.
A quantum vertex algebra is a weak quantum vertex algebra equipped
with a ``rational unitary quantum Yang-Baxter operator'' $\S(x): V\otimes
V\rightarrow V\otimes V\otimes \C((x))$ such that for $u,v\in V$,
(\ref{e1.1}) holds with
$$\sum_{i=1}^{r}v^{(i)}\otimes u^{(i)}\otimes f_{i}(x)
=\S(x)(v\otimes u).$$
The notion of quantum vertex algebra clearly generalizes that of vertex
algebra, vertex superalgebra and vertex color-algebra.

Just as ordinary vertex algebras can be constructed by using 
local sets of (formal) vertex operators (see \cite{li-local}), it was proved in
\cite{li-qva} that weak quantum vertex
algebras can be constructed by using ``$\S$-local'' sets of (formal) vertex operators.
Let $W$ be any vector space. A subset $U$ of
$\Hom (W,W((x)))$, alternatively denoted by $\E(W)$,
is said to be $\S$-local if for any $a(x),b(x)\in U$, there exists 
$$\sum_{i=1}^{r}b^{(i)}(x)\otimes a^{(i)}(x)\otimes f_{i}(x)\in
U\otimes U\otimes \C((x))$$
such that for some nonnegative integer $k$,
\begin{eqnarray}
(x_{1}-x_{2})^{k}a(x_{1})b(x_{2})=(x_{1}-x_{2})^{k}
\sum_{i=1}^{r}f_{i}(x_{2}-x_{1})b^{(i)}(x_{2})a^{(i)}(x_{1}),
\end{eqnarray}
where each $f_{i}(x_{2}-x_{1})$ is
to be expanded in nonnegative powers of $x_{1}$. 
For $a(x),b(x)$ from an $\S$-local set $U$ with the above information,
we define $a(x)_{n}b(x)\in \E(W)$ for $n\in \Z$ as follows:
\begin{eqnarray}
\lefteqn{a(x)_{n}b(x)}\nonumber\\
&=&\Res_{x_{1}}\left(\! (x_{1}-x_{2})^{n}a(x_{1})b(x_{2})-\sum_{i=1}^{r}(-x_{2}+x_{1})^{n}f_{i}(x_{2}-x_{1})
b^{(i)}(x_{2})a^{(i)}(x_{1})\!\right)\!.\ \ \ \
\end{eqnarray}
Then it was proved that for any $\S$-local subset $U$  of $\E(W)$, there exists a unique closed
$\S$-local subspace $\<U\>$ containing $U$ and $1_{W}$ (the identity operator) and that
$\<U\>$ is a weak quantum vertex algebra with $W$ naturally as a module.
Furthermore, it is proved in \cite{li-qva} and this paper
that if $W$ and $U$ are chosen to satisfy certain special properties, 
$\<U\>$ and $W$ are isomorphic $\<U\>$-modules. Consequently,
$W$ has a  weak quantum vertex algebra structure transported from $\<U\>$.
This is analogous to a theorem of \cite{fkrw} and \cite{mp}.
Next, to get a quantum vertex algebra from a weak quantum vertex
algebra, one shall need to construct
a unitary rational quantum Yang-Baxter operator with the desired
property. For this we use an important idea of Etingof and Kazhdan.

In \cite{ek}, Etingof and Kazhdan introduced and studied a notion
of quantum vertex operator algebra in the sense of formal deformation.
By definition, the underlying space of a quantum
vertex operator algebra $V$ is a topologically free $\C[[h]]$-module
and the classical limit $V/hV$ is an ordinary vertex algebra.  Among
the important ingredients of a quantum vertex operator algebra is a
rational quantum Yang-Baxter operator $\S(x)$ (also depending on the
formal parameter $h$).  In that paper, they introduced a notion of
nondegeneracy for a vertex algebra and proved that in constructing a
quantum vertex operator algebra $V$ (in the sense of \cite{ek}), if
$V/hV$ (a vertex algebra over $\C$) is nondegenerate, the quantum
Yang-Baxter relation follows from $\S$-locality and associativity.  
An analogue which was obtained in \cite{li-qva} is that every nondegenerate
weak quantum vertex algebra is a quantum vertex algebra.

Now we have seen that nondegeneracy of (ordinary) vertex
algebras are very important for constructing
quantum vertex operator algebras over $\C[[h]]$ in the sense of
Etingof and Kazhdan while nondegeneracy of nonlocal vertex
algebras are very important for constructing
quantum vertex algebras over $\C$.
The main purpose of this paper is to establish general
results on nondegeneracy of nonlocal vertex algebras, which
can be used to construct interesting examples of quantum vertex algebras.

A nonlocal vertex algebra $V$ is {\em nondegenerate}
if for each positive integer $n$ the linear map
$Z_{n}:V^{\otimes n}\otimes \C((x_{1}))\cdots ((x_{n}))\rightarrow
V((x_{1}))\cdots ((x_{n}))$, defined by
\begin{eqnarray*}
Z_{n}(u_{1}\otimes \cdots \otimes u_{n}\otimes f)=
f Y(u_{1},x_{1})\cdots Y(u_{n},x_{n}){\bf 1},
\end{eqnarray*}
is injective. Let $\hat{\g}$ be an affine Lie algebra and 
let $\ell$ be any complex
number not equal to negative the dual Coxeter number of $\g$. 
Associated to the pair $(\hat{\g},\ell)$ we have 
a vertex operator algebra $V_{\hat{\g}}(\ell,0)$ (cf. \cite{fz}),
whose underlying space is the generalized Verma module or Weyl module of
level $\ell$ for $\hat{\g}$, induced from the $1$-dimensional trivial
$\g$-module $\C$.  It was proved in \cite{ek} that if
$V_{\hat{\g}}(\ell,0)$ is an irreducible $\hat{\g}$-module, then the
vertex operator algebra $V_{\hat{\g}}(\ell,0)$ is
nondegenerate. Later, it was proved in \cite{li-simple} that
for a general vertex algebra $V$, if $V$ is of countable dimension
over $\C$ and if the adjoint module $V$ is irreducible, then $V$ is
nondegenerate. This in particular implies that 
every simple vertex operator algebra in the sense of 
\cite{flm} and \cite{fhl} is nondegenerate. 

In this paper, we first extend the result of \cite{li-simple}, proving that 
if a nonlocal vertex algebra $V$ is countable dimensional over $\C$ and
if the adjoint module $V$ is irreducible, then $V$ is nondegenerate.
Note that for certain rational numbers $\ell$, even though
the vertex operator algebras $V_{\hat{\g}}(\ell,0)$ are not simple, they
play a very important role in many studies.
For the purpose of constructing quantum deformation of 
those nonsimple vertex operator algebras $V_{\hat{\g}}(\ell,0)$
in the sense of \cite{ek}, we certainly hope that they
are also nondegenerate. On the other hand, to apply
the general construction theorem (Theorem \ref{tweak-qva-construction})
of quantum vertex algebras over $\C$,
we often take the underlying space $V$ to be a universal module
for some algebra in a certain sense.
(It seems that interesting nonlocal vertex
algebras are most likely to be reducible, as
most interesting noncommutative associative algebras 
are reducible under the left action and even not simple.)  
Motivated by this, and
we prove that every countable-dimensional nonlocal 
vertex algebra over $\C$, which is of PBW type in a
certain sense, is nondegenerate.  In particular, the vertex operator
algebras $V_{\hat{\g}}(\ell,0)$ associated to affine Lie algebras and
those associated to the Virasoro Lie algebra are nondegenerate.  This
result enables us to construct quantum vertex algebras
associated to certain subalgebras of affine Lie algebras.
For example, for any finite-dimensional Lie algebra $\g$, 
as $d/dt$ is a derivation of the Lie algebra $t^{-1}\g[t^{-1}]$,
it can be extended uniquely to a derivation $d$ of
the universal enveloping algebra $U(t^{-1}\g[t^{-1}])$. 
Then $U(t^{-1}\g[t^{-1}])$ has a nonlocal vertex algebra structure
with $Y(a,x)b=(e^{xd}a)b$ for $a,b\in U(t^{-1}\g[t^{-1}])$
(\cite{bor}, cf. \cite{bk}, \cite{li-g1}).
We prove that this rather trivial nonlocal vertex algebra
actually is a nondegenerate quantum vertex algebra
with a nontrivial quantum Yang-Baxter operator.
It is expected that from Yangians $Y(\g)$, which are closely related to $U(\g[t])$,
we have nondegenerate quantum vertex algebras. This will be addressed in a sequel.

Just as  one needs
associative or Lie algebras of a certain type such as affine Lie
algebras and the Virasoro Lie algebra to construct interesting examples of vertex algebras, 
one needs algebras of a certain type to construct quantum vertex algebras. 
Indeed, there is a family of algebras, called
Zamolodchikov-Faddeev algebras, which are very close to what we need.
Zamolodchikov-Faddeev algebras are associative $C^{*}$-algebras, which 
have appeared in the study of quantum field
theory (see \cite{za}, \cite{fa}, cf. \cite{mr}).  A
Zamolodchikov-Faddeev algebra by definition is associated to an $S$-matrix, or a
quantum Yang-Baxter operator (with two spectral parameters) on a
finite-dimensional vector space $U$.
Motivated by the construction of Zamolodchikov-Faddeev algebras,
{}from any rational unitary quantum
Yang-Baxter operator $\S(x): H\otimes H\rightarrow H\otimes H\otimes
\C[[x]]$ for a finite-dimensional vector space $H$ with $\S(0)=1$
we construct a vector space $V(H,\S)$ and we prove that 
under a certain assumption if $V(H,\S)$ is of PBW type in a certain sense,
$V(H,\S)$ has a unique nondegenerate 
quantum vertex algebra structure with certain properties.
For $(H,\S)$ of a certain special type, we show that
indeed $V(H,\S)$ is of PBW type.

Among our results we also refine the general construction theorem of
\cite{li-qva} for weak quantum vertex algebras.

This paper is organized as follows: In Section 2 we review basic
results on (weak) quantum vertex algebras and we prove certain new
results. In Section 3, we establish certain basic results on
nondegeneracy. In Section 4, we construct a family of quantum vertex
algebras closely related to Zamolodchikov-Faddeev algebras.

\section{Quantum vertex algebras and a general construction}

In this section we first summarize some of the main results of \cite{li-qva} 
on weak quantum vertex algebras and quantum vertex algebras, and we then
give a refinement of a construction theorem obtained in \cite{li-qva}.
We also give an example to illustrate how to use the general construction theorem.

Throughout this paper we use the field of complex numbers $\C$ as the basic field.
Following the tradition we denote by $\C(x_{1},\dots,x_{n})$ 
the field of rational functions.
Set
\begin{eqnarray}
\C_{*}(x_{1},\dots,x_{n})=\{ f/g\;|\; f\in \C[[x_{1},\dots,x_{n}]],\; 
0\ne g\in \C[x_{1},\dots,x_{n}]\},
\end{eqnarray}
a subring of the fraction field of $\C[[x_{1},\dots,x_{n}]]$.
So $\C_{*}(x_{1},\dots,x_{n})$ contains 
the field of rational functions $\C(x_{1},\dots,x_{n})$ as a subring.
Denote by $\iota_{x_{1},x_{2}}$ the natural embedding of
$\C_{*}(x_{1},x_{2})$ into the field $\C((x_{1}))((x_{2}))$.

Throughout this paper, nonlocal vertex algebras
are synonymous to $G_{1}$-vertex algebras studied in \cite{li-g1}
and they are also essentially field algebras studied 
in \cite{bk} (cf. \cite{kac}).

\bd{dgva}
{\em A {\em nonlocal vertex algebra} is a vector space $V$ equipped
with a linear map 
\begin{eqnarray}
Y(\cdot,x): V &\rightarrow & 
\Hom (V,V((x)))\subset (\End V)[[x,x^{-1}]]\nonumber\\
v&\mapsto& Y(v,x)=\sum_{n\in \Z}v_{n}x^{-n-1}\;\;\; (v_{n}\in \End V)
\end{eqnarray}
and equipped with a distinguished vector ${\bf 1}$ such that
for $v\in V$
\begin{eqnarray}
& &Y({\bf 1},x)v=v,\\
& &Y(v,x){\bf 1}\in V[[x]]\;\;\mbox{ and }\;\; \lim_{x\rightarrow
0}Y(v,x){\bf 1}=v
\end{eqnarray}
and such that for $u,v,w\in V$, there exists a nonnegative integer
$l$ such that
\begin{eqnarray}
(x_{0}+x_{2})^{l}Y(u,x_{0}+x_{2})Y(v,x_{2})w=
(x_{0}+x_{2})^{l}Y(Y(u,x_{0})v,x_{2})w.
\end{eqnarray}}
\ed

\bd{dmodule}
{\em A $V$-{\em module} is a vector space $W$ equipped with a linear map
\begin{eqnarray}
Y_{W}(\cdot,x): V &\rightarrow & 
\Hom (W,W((x)))\subset (\End W)[[x,x^{-1}]]\nonumber\\
v&\mapsto& Y_{W}(v,x)=\sum_{n\in \Z}v_{n}x^{-n-1}\;\;\; (v_{n}\in \End W)
\end{eqnarray}
such that all the following axioms hold: 
\begin{eqnarray}
Y_{W}({\bf 1},x)=1_{W}\;\;\;
\mbox{(where $1_{W}$ is the identity operator on $W$)},
\end{eqnarray}
and for any $u, v\in V,\; w\in W$, there exists $l\in \N$
such that
\begin{eqnarray}\label{emoduleweakassoc}
(x_{0}+x_{2})^{l}Y_{W}(u,x_{0}+x_{2})Y_{W}(v,x_{2})w
=(x_{0}+x_{2})^{l}Y_{W}(Y(u,x_{0})v,x_{2})w.
\end{eqnarray}}
\ed

\bd{dweak-qva} {\em A {\em weak quantum vertex algebra} is a
vector space $V$ equipped with a vector ${\bf 1}$ and a linear map
$Y: V\rightarrow \Hom (V,V((x)))$ such that for $v\in V$,
\begin{eqnarray}
& &Y({\bf 1},x)v=v,\\
& &Y(v,x){\bf 1}\in V[[x]]\;\;\mbox{ and }\; \lim_{x\rightarrow 0}Y(v,x){\bf 1}=v
\end{eqnarray}
and such that for any $u,v\in V$, there exists
$\sum_{i=1}^{r}v^{(i)}\otimes u^{(i)}\otimes f_{i}(x)\in V\otimes V\otimes \C((x))$ such that
\begin{eqnarray}\label{es-jacobi}
& &x_{0}^{-1}\delta\left(\frac{x_{1}-x_{2}}{x_{0}}\right)Y(u,x_{1})Y(v,x_{2})
-x_{0}^{-1}\delta\left(\frac{x_{2}-x_{1}}{-x_{0}}\right)
\sum_{i=1}^{r}f_{i}(-x_{0})Y(v^{(i)},x_{2})Y(u^{(i)},x_{1})\nonumber\\
& &=x_{2}^{-1}\delta\left(\frac{x_{1}-x_{0}}{x_{2}}\right)Y(Y(u,x_{0})v,x_{2}).
\end{eqnarray}}
\ed

\br{requivalent}
{\em One can show that (\ref{es-jacobi}) is equivalent to 
the {\em $\S$-locality}: For $u,v\in V$, there exists
$\sum_{i=1}^{r}v^{(i)}\otimes u^{(i)}\otimes f_{i}(x)\in V\otimes V\otimes \C((x))$
such that 
\begin{eqnarray}
(x_{1}-x_{2})^{k}Y(u,x_{1})Y(v,x_{2})
=(x_{1}-x_{2})^{k}\sum_{i=1}^{r}\iota_{x_{2},x_{1}}(f_{i}(x_{2}-x_{1}))Y(v^{(i)},x_{2})Y(u^{(i)},x_{1})
\end{eqnarray}
for some nonnegative integer $k$, 
and the {\em weak associativity}: For any $u,v,w\in V$, there exists a nonnegative integer $l$ such that
\begin{eqnarray}
(x_{0}+x_{2})^{l}Y(u,x_{0}+x_{2})Y(v,x_{2})w=(x_{0}+x_{2})^{l}Y(Y(u,x_{0})v,x_{2})w.
\end{eqnarray}
Therefore a weak quantum vertex algebra is exactly
a nonlocal vertex algebra $V$ which satisfies $\S$-locality.}
\er

\br{rwqva-module}
{\em Let $V$ be a weak quantum vertex algebra and 
let $(W,Y_{W})$ be a module for $V$ 
viewed as a nonlocal vertex algebra. {}From \cite{li-qva}, 
whenever (\ref{es-jacobi}) holds, we have
\begin{eqnarray}\label{es-jacobi-module}
& &x_{0}^{-1}\delta\left(\frac{x_{1}-x_{2}}{x_{0}}\right)Y_{W}(u,x_{1})Y_{W}(v,x_{2})\nonumber\\
& &\hspace{1cm} \ \ \ -x_{0}^{-1}\delta\left(\frac{x_{2}-x_{1}}{-x_{0}}\right)
\sum_{i=1}^{r}f_{i}(-x_{0})Y_{W}(v^{(i)},x_{2})Y_{W}(u^{(i)},x_{1})\nonumber\\
& &=x_{2}^{-1}\delta\left(\frac{x_{1}-x_{0}}{x_{2}}\right)Y_{W}(Y(u,x_{0})v,x_{2}).
\end{eqnarray}
In view of this, we simply define a $V$-module to be a module for $V$ 
viewed as a nonlocal vertex algebra. }
\er

Let $W$ be any vector space. Set
$$\E(W)=\Hom (W,W((x))),$$
which in particular contains the identity operator $1_{W}$ on $W$.
A subset $U$ of $\E(W)$ is said to be {\em $\S$-local}
if for any $a(x),b(x)\in U$, there exist
$f_{i}(x)\in \C((x)),\; b^{(i)}(x),a^{(i)}(x)\in U$ for $i=1,\dots,r$, 
and a nonnegative integer $k$ such that 
\begin{eqnarray}\label{eS-local-ab}
(x_{1}-x_{2})^{k}a(x_{1})b(x_{2})
=(x_{1}-x_{2})^{k}\sum_{i=1}^{r}
\iota_{x_{2},x_{1}}(f_{i}(x_{2}-x_{1}))b^{(i)}(x_{2})a^{(i)}(x_{1}).
\end{eqnarray}
Notice that the above relation implies 
\begin{eqnarray}
(x_{1}-x_{2})^{k}a(x_{1})b(x_{2})\in \Hom (W,W((x_{1},x_{2}))).
\end{eqnarray}
Let $U$ be any $\S$-local subset of $\E(W)$. For $a(x),b(x)\in U$ and for $n\in \Z$
we define $a(x)_{n}b(x)\in \E(W)$ in terms of the generating function
\begin{eqnarray}
Y_{\E}(a(x),x_{0})b(x)=\sum_{n\in \Z}(a(x)_{n}b(x))x_{0}^{-n-1}
\end{eqnarray}
 by
\begin{eqnarray}
Y_{\E}(a(x),x_{0})b(x)=\Res_{x_{1}}x^{-1}\delta\left(\frac{x_{1}-x_{0}}{x}\right)x_{0}^{-k}
\left((x_{1}-x)^{k}a(x_{1})b(x)\right),
\end{eqnarray}
where $k$ is any nonnegative integer such that 
$(x_{1}-x_{2})^{k}a(x_{1})b(x_{2})\in \Hom (W,W((x_{1},x_{2})))$. 
Assuming the $\S$-locality relation (\ref{eS-local-ab}) we have
\begin{eqnarray}
& &Y_{\E}(a(x),x_{0})b(x)
=\Res_{x_{1}}x_{0}^{-1}\delta\left(\frac{x_{1}-x}{x_{0}}\right)a(x_{1})b(x)\nonumber\\
& &\hspace{2cm}-\Res_{x_{1}}x_{0}^{-1}\delta\left(\frac{x-x_{1}}{-x_{0}}\right)\sum_{i=1}^{r}
\iota_{x,x_{1}}(f_{i}(x-x_{1}))b^{(i)}(x)a^{(i)}(x_{1}),
\end{eqnarray}
or equivalently, for $n\in \Z$,
\begin{eqnarray}
\lefteqn{a(x)_{n}b(x)}\nonumber\\
&=&\Res_{x_{1}}\left(\! (x_{1}-x_{2})^{n}a(x_{1})b(x_{2})-\sum_{i=1}^{r}(-x_{2}+x_{1})^{n}f_{i}(x_{2}-x_{1})
b^{(i)}(x_{2})a^{(i)}(x_{1})\!\right)\!.\ \ \ \
\end{eqnarray}
An $\S$-local subspace $U$ is said to be {\em closed} if 
$a(x)_{n}b(x)\in U$ for $a(x),b(x)\in U,\; n\in \Z$.
The following is one of the results obtained in \cite{li-qva}:

\bt{tX}
Let $W$ be any vector space and let $U$ be any $\S$-local subset of $\E(W)$.
Then there exists a (unique) smallest
closed $\S$-local subspace $\<U\>$, containing $U$ and $1_{W}$. Furthermore,
$(\<U\>,Y_{\E},1_{W})$ is a weak quantum vertex algebra and $(W,Y_{W})$ carries the structure of
a faithful $V$-module where $Y_{W}(\alpha(x),x_{0})=\alpha(x_{0})$ for $\alpha(x)\in V$.
\et

Let $V$ be a nonlocal vertex algebra. A subset $U$ of $V$
is said to be {\em $S$-local}
if $\{ Y(u,x)\;|\; u\in U\}$ is an $\S$-local subset of $\E(V)$.
Note that a subset of an $\S$-local subset is not necessarily $\S$-local.
Nevertheless we have:

\bl{lS-local}
Let $V$ be a nonlocal vertex algebra and let $U$ be an $\S$-local subspace.
Denote by $U^{(2)}$ the linear span of $U$ and the vectors 
$a_{m}b$ for $a,b\in U,\; m\in \Z$.
Then $U^{(2)}$ is $\S$-local. Furthermore, the subalgebra $\<U\>$
generated by $U$ is a weak quantum vertex algebra.
\el

\begin{proof} By assumption $\overline{U}=\{ Y(v,x)\;|\; v\in V\}$ 
is an $\S$-local subset of $\E(V)$.
It follows from \cite{li-qva} (Lemmas 5.7 and 5.8) that
$\overline{U}^{(2)}$ is $\S$-local.
For $u,v\in U$, by Proposition 6.6 of \cite{li-qva}
(noting that $\S$-locality implies compatibility), we have
$$Y(u_{m}v,x)=Y(u,x)_{m}Y(v,x).$$
It follows that $\overline{U^{(2)}}=\overline{U}^{(2)}$.
Thus $\overline{U^{(2)}}$ is $\S$-local.
By definition $U^{(2)}$ is $\S$-local.

For the second assertion, let $K$ be a maximal $\S$-local subspace
of $\<U\>$, containing $U$. Clearly, $K+\C {\bf 1}$ is $\S$-local.
As $K$ is maximal we have ${\bf 1}\in K$.
By the first assertion, $K^{(2)}$ is $\S$-local.
Again, as $K$ is maximal we have $K^{(2)}=K$. Thus $K$ is a
subalgebra of $\<U\>$, containing ${\bf 1}$ and $U$. 
Therefore $\<U\>=K$ and it is a weak quantum vertex algebra.
\end{proof}

We shall need the following generalization of a result 
in \cite{li-form} and \cite{li-qva}:

\bp{pvacuum}
Let $V$ be a weak quantum vertex algebra, let $W$ be a $V$-module and
let $e\in W$ be such that
$Y(u,x)e\in W[[x]]$ for $u\in U$, where $U$ is an $\S$-local generating subspace of $V$.
Then $Y(v,x)e\in W[[x]]$ for all $v\in V$. Furthermore, the linear map from
$V$ to $W$, sending $v$ to $v_{-1}e$, is a $V$-module homomorphism.
\ep

\begin{proof} {}From \cite{li-qva} (cf. \cite{li-form}),
the second assertion is true if the first assertion is proved.
Let $F$ be the collection of $\S$-local subspaces $L$ containing $U$ and 
satisfying the condition $Y(a,x)e\in W[[x]]$ for $a\in L$.
By the assumption, $U\in F$. Let $K$ be a maximal element of $F$.
We are going to prove $K=V$. By assumption, $U\subset K$.
Clearly $K+\C{\bf 1}\in F$.
It follows that ${\bf 1}\in K$. As $U$ generates $V$ as a weak quantum vertex algebra, 
it suffices to prove that $K$ is closed.

For any $a,b\in K$, there exist $a^{(i)}, b^{(i)}\in K,\; f_{i}(x)\in \C((x))$ 
for $i=1,\dots,r$ such that
\begin{eqnarray}
& &x_{2}^{-1}\delta\left(\frac{x_{1}-x_{0}}{x_{2}}\right)Y(Y(a,x_{0})b,x_{2})e
-x_{0}^{-1}\delta\left(\frac{x_{1}-x_{2}}{x_{0}}\right)Y(a,x_{1})Y(b,x_{2})e
\nonumber\\
&=&-x_{0}^{-1}\delta\left(\frac{x_{2}-x_{1}}{-x_{0}}\right)
\sum_{i=1}^{r}\iota_{x_{2},x_{1}}(f_{i}(x_{2}-x_{1}))Y(b^{(i)},x_{2})Y(a^{(i)},x_{1})e.
\end{eqnarray}
As $Y(a^{(i)},x_{1})e\in W[[x_{1}]]$ for all the $i$'s, taking $\Res_{x_{1}}$ we get
$$Y(Y(a,x_{0})b,x_{2})e
=\Res_{x_{1}}x_{0}^{-1}\delta\left(\frac{x_{1}-x_{2}}{x_{0}}\right)Y(a,x_{1})Y(b,x_{2})e
=Y(a,x_{0}+x_{2})Y(b,x_{2})e.$$
Since $Y(b,x_{2})e\in W[[x_{2}]]$, the right-hand side lies in
$W((x_{0}))[[x_{2}]]$. Consequently, we have
$$Y(Y(a,x_{0})b,x_{2})e\in W((x_{2}))((x_{0}))\cap W((x_{0}))[[x_{2}]]
\subset W[[x_{0},x_{2}]][x_{0}^{-1}].$$
Thus $Y(a_{m}b,x)e\in W[[x]]$ for $m\in \Z$. 
On the other hand, by Lemma \ref{lS-local}, $K^{(2)}$ is $\S$-local. 
Therefore $K^{(2)}\in F$.
Since $K$ is maximal and $K\subset K^{(2)}$, we have $K^{(2)}=K$. 
This proves that $K$ is closed and concludes the proof.
\end{proof}

The following is our general construction theorem
(cf. \cite{li-qva}, \cite{fkrw}, \cite{mp}:

\bt{tweak-qva-construction}
Let $V$ be a vector space, $U$ a subset, ${\bf 1}$ a vector of $V$,
and $Y_{0}$ a map from $U$ to
$\E(V)\;(=\Hom (V,V((x))))$. Assume all the following conditions hold:
\begin{eqnarray}\label{eqva-creation}
& &Y_{0}(u,x){\bf 1}\in V[[x]]\;\;\;\mbox{ and }\;\;
\lim_{x\rightarrow 0}Y_{0}(u,x){\bf 1}=u\;\;\;\mbox{ for }u\in U,
\end{eqnarray}
$Y_{0}(U)=\{ Y_{0}(u,x)\;|\; u\in U\}\subset \E(V)$ is $\S$-local, and
\begin{eqnarray}
V={\rm span}\{ u^{(1)}_{n_{1}}\cdots u^{(r)}_{n_{r}}{\bf 1}\;|\;
r\ge 0,\; u^{(i)}\in U,\; n_{i}\in \Z\},
\end{eqnarray}
where $Y_{0}(u,x)=\sum_{n\in \Z}u_{n}x^{-n-1}$ for $u\in U$.
Furthermore, we assume that there exists a linear map
$f$ from $V$ to $\<Y_{0}(U)\>$ such that $f({\bf 1})=1_{V}$ and
\begin{eqnarray}
f(u_{n}v)=u(x)_{n}f(v)
\;\;\;\mbox{ for }u\in U,\; n\in \Z,\; v\in V.
\end{eqnarray}
Then $Y_{0}$ can be extended uniquely to a linear map $Y$ from $V$ to $\E(V)$ 
such that $(V,Y,{\bf 1})$ carries the structure of 
a weak quantum vertex algebra.
\et

\begin{proof} We shall just prove the existence 
as the uniqueness is clear.
By assumption, $Y_{0}(U)$ is an $\S$-local subset of $\E(V)$.
By Theorem \ref{tX}, $Y_{0}(U)$ generates 
a weak quantum vertex algebra $\<Y_{0}(U)\>$
inside $\E(V)$ with $V$ as a module where
$Y_{V}(\alpha(x),x_{0})=\alpha(x_{0})$ for $\alpha(x)\in \<Y_{0}(U)\>$. 
With the assumption (\ref{eqva-creation}), by Proposition \ref{pvacuum}, 
the linear map $g$ from $\<Y_{0}(U)\>$ to $V$, 
sending $\alpha(x)$ to $\Res_{x}x^{-1}\alpha(x){\bf 1}$,
is a $\<Y_{0}(U)\>$-module homomorphism.
Clearly, $g$ sends $1_{V}$ to ${\bf 1}$.
Consequently, $f$ and $g$ are inverse each other.
Now $V$ has a weak quantum vertex algebra structure $Y$ transported from
$\<Y_{0}(U)\>$ through the linear isomorphism $f$.
For $u\in U$, as
$$f(u)=f(u_{-1}{\bf 1})=u(x)_{-1}1_{V}=u(x),$$
we have
$$Y(u,x_{0})=gY_{\E}(u(x),x_{0})f=Y_{V}(u(x),x_{0})=u(x_{0}).$$
Indeed, $Y$ extends the linear map $Y_{0}$.
This completes the proof.
\end{proof}

\br{rrefine}
{\em In \cite{li-qva}, a general construction theorem for nonlocal vertex algebras
was given, which is
an analogue of a construction theorem for vertex algebras of Frenkel-Kac-Ratiu-Wang \cite{fkrw}
and Meurman-Primc \cite{mp}.
In \cite{li-qva}, $Y_{0}(U)$ was assumed to be quasi compatible (weaker than $\S$-locality)
and in addition it was assumed that
there exists a linear operator $L(-1)$ on $V$ such that
\begin{eqnarray*}
& &L(-1){\bf 1}=0,\\
& &[L(-1),Y_{0}(u,x)]=\frac{d}{dx}Y_{0}(u,x)\;\;\;\mbox{ for }u\in U.
\end{eqnarray*}}
\er

The following result, obtained in \cite{li-qva}, 
is an analogue of Theorem 5.7.6 of \cite{ll}:

\bt{tZ} Let $V$ be a weak quantum vertex algebra, let $U$ be a generating
subspace of $V$, let $W$ be a vector space, and let $Y_{W}^{o}$ be a
linear map from $U$ to $\E(W)$. Suppose that $Y_{W}^{o}(U)$ is
$\S$-local and that there exists a linear map $f$ from
$V$ to $\<Y_{W}^{o}(U)\>$ such that $f({\bf 1})=1_{W}$ and
\begin{eqnarray}
f(u_{n}v)=u(x)_{n}f(v)\;\;\;\mbox{ for }
u\in U,\; n\in \Z,\; v\in V.
\end{eqnarray}
Then $Y_{W}^{o}$ can be uniquely extended to a linear map 
$Y_{W}$ from $V$ to $\E(W)$ such that $(W,Y_{W})$ 
carries the structure of a $V$-module.
\et

In practice we shall use the following result of \cite{li-qva}
as a companion of Theorems \ref{tweak-qva-construction} and \ref{tZ}:

\bp{pY}
Let $W$ be a vector space, let $V$ be a closed $\S$-local subspace 
of $\E(W)$, and let
\begin{eqnarray*}
& &u(x), v(x), u^{(1)}(x), v^{(1)}(x),\dots, u^{(r)}(x), v^{(r)}(x), c^{0}(x),\dots, c^{s}(x) \in V,\\
& &f_{1}(x),\dots, f_{r}(x)\in \C((x)).
\end{eqnarray*}
Suppose that
\begin{eqnarray}\label{e-cross-bracket-W}
& &(x_{1}-x_{2})^{n}u(x_{1})v(x_{2})
-(-x_{2}+x_{1})^{n}\sum_{i=1}^{r}f_{i}(x_{2}-x_{1})v^{(i)}(x_{2})u^{(i)}(x_{1})
\nonumber\\
&=&
\sum_{j=0}^{s}c^{j}(x_{2})\frac{1}{j!}\left(\frac{\partial}{\partial x_{2}}\right)^{j}
x_{2}^{-1}\delta\left(\frac{x_{1}}{x_{2}}\right)
\end{eqnarray}
for some integer $n$. Then 
\begin{eqnarray}\label{e-cross-bracket-algebra-case}
& &(x_{1}-x_{2})^{n}Y_{\E}(u(x),x_{1})Y_{\E}(v(x),x_{2})\nonumber\\
& &\ \ \ \ \ \ 
-(-x_{2}+x_{1})^{n}\sum_{i=1}^{r}f_{i}(x_{2}-x_{1})Y_{\E}(v^{(i)}(x),x_{2})Y_{\E}(u^{(i)}(x),x_{1})
\nonumber\\
&=&\sum_{j=0}^{s}Y_{\E}(c^{j}(x),x_{2})\frac{1}{j!}\left(\frac{\partial}{\partial x_{2}}\right)^{j}
x_{2}^{-1}\delta\left(\frac{x_{1}}{x_{2}}\right).
\end{eqnarray}
\ep

\bd{dqva} {\em A
{\em quantum vertex algebra} is a weak quantum vertex algebra $V$
equipped with a unitary rational quantum Yang-Baxter operator $\S$, 
which by definition is a linear map from $V\otimes V$ to 
$V\otimes V\otimes \C((x))$ satisfying the unitarity relation
\begin{eqnarray}
\S_{21}(-x)\S(x)=1
\end{eqnarray}
and the following quantum Yang-Baxter equation
\begin{eqnarray}
\S_{12}(x)\S_{13}(x+z)\S_{23}(z)=\S_{23}(z)\S_{13}(x+z)\S_{12}(x),
\end{eqnarray}
such that for $u,v\in V$, (\ref{es-jacobi}) holds with
$\sum_{i=1}^{r}v^{(i)}\otimes u^{(i)}\otimes f_{i}(x)=\S(x)(v\otimes u)$.} 
\ed

The following result was proved in \cite{li-qva}
(cf. \cite{ek}):

\bp{pweak-nonweak}
Let $V$ be a weak quantum vertex algebra. 
Assume that the linear map 
$\pi_{3}: V\otimes V\otimes V\rightarrow 
\Hom (V,V((x_{1}))((x_{2}))((x_{3})))$, defined by
\begin{eqnarray}
\pi_{3}(u\otimes v\otimes w)=Y(u,x_{1})Y(v,x_{2})Y(w,x_{3})
\end{eqnarray}
for $u,v,w\in V$, is injective. 
For $u,v\in V$, we define 
$$\S(x)(v\otimes
u)=\sum_{i=1}^{r}a_{i}\otimes b_{i}\otimes f_{i}(x)\in V\otimes
V\otimes \C((x))$$
if
\begin{eqnarray}
(x_{1}-x_{2})^{k}Y(u,x_{1})Y(v,x_{2})
=(x_{1}-x_{2})^{k}\sum_{i=1}^{r}\iota_{x_{2},x_{1}}(f_{i}(x_{2}-x_{1}))
Y(a_{i},x_{2})Y(b_{i},x_{1})
\end{eqnarray}
for some nonnegative integer $k$.
Then $\S(x)$ is well defined and it is a rational unitary quantum Yang-Baxter
operator on $V$ and furthermore, $(V,\S)$ is a quantum vertex algebra.
\ep

Closely related the linear map $\pi_{3}$ is the following notion 
of nondegeneracy due to Etingof and Kazhdan (see \cite{ek}):

\bd{dek-voa}
{\em A nonlocal vertex algebra $V$ is said to be {\em nondegenerate}
if for each positive integer $n$ the linear map
$Z_{n}: V^{\otimes n}\otimes \C((x_{1}))\cdots ((x_{n}))\rightarrow
V((x_{1}))\cdots ((x_{n}))$, defined by
\begin{eqnarray}
Z_{n}(u_{1}\otimes \cdots \otimes u_{n}\otimes f)=
f Y(u_{1},x_{1})\cdots Y(u_{n},x_{n}){\bf 1}
\end{eqnarray}
is injective.}
\ed

Clearly, $\pi_{3}$ is injective whenever $Z_{3}$ is injective. Therefore we have:

\bt{tweak-strongqva}
Every nondegenerate weak quantum vertex algebra is a quantum vertex algebra 
where the quantum Yang-Baxter operator $\S(x)$ is uniquely determined.
\et

In the following we illustrate how 
to apply the general construction 
(Theorems \ref{tX} and \ref{tweak-qva-construction} and Proposition \ref{pY})
to construct weak quantum vertex algebras.

\bd{dQLie-module}
{\em Let $H$ be a vector space 
equipped with a bilinear form $\<\cdot,\cdot\>$ 
and let $\S (x): H\otimes H\rightarrow H\otimes H\otimes \C((x))$ be a 
linear map. An {\em $(H,\S)$-module} is 
a module $W$ for the (free) tensor algebra $T(H\otimes \C[t,t^{-1}])$ 
such that
for any $a\in H,\; w\in W$,
\begin{eqnarray}\label{e-truncation-qLie}
a(m)w=0\;\;\;\mbox{ for $m$ sufficiently large},
\end{eqnarray}
where we use $a(m)$ for the operator corresponding to $a\otimes t^{m}$,  and such that
\begin{eqnarray}\label{emain-relation-qLie}
a(x_{1})b(x_{2})w
-\sum_{i=1}^{r}\iota_{x_{2},x_{1}}(f_{i}(x_{2}-x_{1}))b^{(i)}(x_{2})a^{(i)}(x_{1})w
=x_{2}^{-1}\delta\left(\frac{x_{1}}{x_{2}}\right)\<a,b\>w
\end{eqnarray}
for $a,b\in H,\; w\in W$, where 
$\S(x)(b\otimes a)=\sum_{i=1}^{r}b^{(i)}\otimes a^{(i)}\otimes f_{i}(x)$. 
(Notice that the condition (\ref{e-truncation-qLie}) guarantees that
 for any $m,n\in \Z$,
the coefficient of $x_{1}^{m}x_{2}^{n}$ in the left side is a finite sum in $W$.)}
\ed

Let $W$ be an $(H,\S)$-module.
A vector $e\in W$ is called a {\em vacuum-like vector}  if 
\begin{eqnarray}
a(m)e=0\;\;\;\mbox{ for all }a\in H,\; m\ge 0.
\end{eqnarray}
We call $W$ a {\em vacuum $(H,\S)$-module} if
$W$ is cyclic on a vacuum-like vector.

\bp{pqLie-wqva}
Let $H$ be a vector space equipped with a bilinear form $\<\cdot,\cdot\>$ and 
let $\S: H\otimes H\rightarrow H\otimes H\otimes \C((x))$ be a linear map.
Suppose that $V(H,\S)$ is a universal vacuum $(H,\S)$-module with
a vacuum vector ${\bf 1}$ as its generator such that
the linear map from $H$ to $V(H,\S)$, sending
$a$ to $a(-1){\bf 1}$ for $a\in H$, is injective.
Then there exists a unique
weak quantum vertex algebra structure on $V(H,\S)$ with ${\bf 1}$ 
as the vacuum vector such that
\begin{eqnarray}
Y(a(-1){\bf 1},x)=a(x)\;\;\;\mbox{ for }a\in H.
\end{eqnarray}
Furthermore, on any $(H,\S)$-module $W$ there exists a unique $V(H,\S)$-module
structure $Y_{W}$ such that
\begin{eqnarray}
Y_{W}(a(-1){\bf 1},x)=a(x)\;\;\;\mbox{ for }a\in H.
\end{eqnarray}
\ep

\begin{proof} We are going to apply Theorem \ref{tweak-qva-construction}.
Take $V=V(H,\S)$ and $U=\{ a(-1){\bf 1}\;|\; a\in H\}$.
Define $Y_{0}(a(-1){\bf 1},x)=a(x)$ for $a\in H$. 
The relation (\ref{emain-relation-qLie}) implies that
$\{ Y_{0}(a(-1){\bf 1},x)=a(x)\;|\; a\in H\}$ is $\S$-local.
Then all the assumptions except 
the last one (about the existence of a certain linear map $f$) hold.

For $a,b\in H$, with the relation (\ref{emain-relation-qLie}) 
by Proposition \ref{pY} we have
\begin{eqnarray*}
& &Y_{\E}(a(x),x_{1})Y_{\E}(b(x),x_{2})
-\sum_{i=1}^{r}\iota_{x_{2},x_{1}}(f_{i}(x_{2}-x_{1}))Y_{\E}(b^{(i)}(x),x_{2})Y_{\E}(a^{(i)}(x),x_{1})
\nonumber\\
&=&x_{2}^{-1}\delta\left(\frac{x_{1}}{x_{2}}\right)\<a,b\>.
\end{eqnarray*}
Thus $\<Y_{0}(U)\>$ is an $(H,\S)$-module with $a(x_{0})$ acting as $Y_{\E}(a(x),x_{0})$ for $a\in H$.
Furthermore, $\<Y_{0}(U)\>$ is a vacuum module with the vacuum vector $1_{V}$ as its generator.
As $V(H,\S)$ is a universal vacuum $(H,\S)$-module,
there exists an
$(H,\S)$-homomorphism $f$ from $V(H,\S)$ to $\<Y_{0}(U)\>$,
sending ${\bf 1}$ to $1_{V}$. 
Now it follows from Theorem \ref{tweak-qva-construction}. The last assertion follows from
Theorem \ref{tZ}.
\end{proof}

\br{rgeneral-case}
{\em If the linear map $\pi_{3}$ defined in Proposition \ref{pweak-nonweak}
for $V(H,\S)$ is injective, it follows from Proposition
\ref{pweak-nonweak} that $\S$ is necessarily a unitary rational 
quantum Yang-Baxter operator on $H$. In Section 4, we shall study 
this again after we establish certain results on nondegeneracy in Section 3.}
\er

\br{rusual-ZFalgebra}
{\em The notion of $(H,\S)$-module was motivated by the following notion 
of Zamolodchikov-Faddeev algebra (see [Z], [F]). Let $U$ be a finite-dimensional
vector space over $\C$ and let $\S$ be
a quantum Yang-Baxter operator on $U$, which by definition is a 
linear map 
$$\S(x_{1},x_{2}): U\otimes U\rightarrow U\otimes U\otimes \C_{*}(x_{1},x_{2}),$$
satisfying the quantum Yang-Baxter equation
\begin{eqnarray}
\S_{12}(x_{1},x_{2})\S_{13}(x_{1},x_{3})\S_{23}(x_{2},x_{3})
=\S_{23}(x_{2},x_{3})\S_{13}(x_{1},x_{3})\S_{12}(x_{1},x_{2}).
\end{eqnarray}
The Zamolodchikov-Faddeev (ZF) algebra $\A (U,\S)$ was defined 
to be an associative ($C^{*}$-)algebra with identity ${\bf 1}$ with generators
$\{ a^{*\alpha}(\chi),a_{\alpha}(\chi)\;|\; \chi\in \R\}$ 
satisfying the following relations:
\begin{eqnarray}
& &a_{\alpha_{1}}(\chi_{1})a_{\alpha_{2}}(\chi_{2})
-\S_{\alpha_{2}\alpha_{1}}^{\beta_{2}\beta_{1}}(\chi_{2},\chi_{1})
a_{\beta_{2}}(\chi_{2})a_{\beta_{1}}(\chi_{1})=0,\\
& &a^{*\alpha_{1}}(\chi_{1})a^{*\alpha_{2}}(\chi_{2})
-\S^{\alpha_{2}\alpha_{1}}_{\beta_{2}\beta_{1}}(\chi_{2},\chi_{1})
a^{*\beta_{2}}(\chi_{2})a^{*\beta_{1}}(\chi_{1})=0,\\
& &a_{\alpha_{1}}(\chi_{1})a^{*\alpha_{2}}(\chi_{2})
-a^{*\beta_{2}}(\chi_{2})\S_{\alpha_{1}\beta_{2}}^{\beta_{1}\alpha_{2}}(\chi_{1},\chi_{2})
a_{\beta_{1}}(\chi_{1})=2\pi \delta_{\alpha_{1}\alpha_{2}}\delta(\chi_{1}-\chi_{2}){\bf 1}.
\ \ \ \
\end{eqnarray}
We will not use this notion in this paper.}
\er

\section{Nondegeneracy of nonlocal vertex algebras}
Having known that
every nondegenerate weak quantum vertex algebra is a 
quantum vertex algebra, in this section we determine 
when a nonlocal vertex algebra is nondegenerate.
We prove that any irreducible nonlocal vertex algebra of countable dimension 
over $\C$ is nondegenerate, generalizing an earlier result of \cite{li-simple}
{}from (ordinary) vertex algebras to nonlocal vertex algebras.
We also prove that any nonlocal vertex algebra with a basis of PBW type is 
nondegenerate. Then we apply these results to establish the nondegeneracy of
better known vertex (operator) algebras. We also construct some 
nondegenerate quantum vertex algebras associated to certain Lie algebras.

Recall that a nonlocal vertex algebra $V$ is said to be nondegenerate
if for each positive integer $n$ the linear map
$Z_{n}:V^{\otimes n}\otimes \C((x_{1}))\cdots ((x_{n}))\rightarrow
V((x_{1}))\cdots ((x_{n}))$, defined by
\begin{eqnarray*}
Z_{n}(u_{1}\otimes \cdots \otimes u_{n}\otimes f)=
f Y(u_{1},x_{1})\cdots Y(u_{n},x_{n}){\bf 1}
\end{eqnarray*}
is injective.

Note that $\C((x_{1}))\cdots ((x_{n}))$ is a field and
$V^{\otimes n}\otimes \C((x_{1}))\cdots ((x_{n}))$ and $V((x_{1}))\cdots ((x_{n}))$ 
are vector spaces over $\C((x_{1}))\cdots ((x_{n}))$. Clearly,
the 
$\C$-linear maps  $Z_{n}$ actually are $\C((x_{1}))\cdots ((x_{n}))$-linear maps.

First we have the following result:

\bl{lnon-tensor}
Let $U$ and $V$ be nondegenerate nonlocal vertex algebras. Then
the tensor product nonlocal vertex algebra $U\otimes V$ is nondegenerate.
\el

\begin{proof} For any positive integer $n$, we have
\begin{eqnarray*}
& &(U\otimes V)^{\otimes n}\otimes \C((x_{1}))\cdots ((x_{n}))\\
&=&(U^{\otimes n}\otimes \C((x_{1}))\cdots ((x_{n})))\otimes_{\C((x_{1}))\cdots ((x_{n}))}(V^{\otimes n}\otimes \C((x_{1}))\cdots ((x_{n})))
\end{eqnarray*}
and the image of $Z_{n}^{U\otimes V}$ (the $Z_{n}$ map $U\otimes V$) is contained in 
\begin{eqnarray*}
U((x_{1}))\cdots ((x_{n}))\otimes_{ \C((x_{1}))\cdots ((x_{n}))}V((x_{1}))\cdots ((x_{n})),
\end{eqnarray*}
a subset of $(U\otimes V)((x_{1}))\cdots ((x_{n}))$.
Viewing $Z_{n}^{U}, Z_{n}^{V}$ and $Z_{n}^{U\otimes V}$ 
as $\C((x_{1}))\cdots ((x_{n}))$-linear maps, we have
$$Z_{n}^{U\otimes V}=Z_{n}^{U}\otimes_{ \C((x_{1}))\cdots ((x_{n}))} Z_{n}^{V}.$$
As both $Z_{n}^{U}$ and $Z_{n}^{V}$ are injective, so is $Z_{n}^{U\otimes V}$.
\end{proof}

Immediately from definition we have:

\bl{lnon-subva}
Let $V$ be a nondegenerate nonlocal vertex algebra.
Then any subalgebra of $V$ is nondegenerate. 
\el

On the other hand we have:

\bl{lnon-kernel}
Let $V$ be a nonlocal vertex algebra and let $\D$ be the linear
operator defined by $\D(v)=v_{-2}{\bf 1}$ for $v\in V$.
If $\dim \ker \D\ge 2$, then $V$ is degenerate. 
In particular, any associative algebra of dimension greater than $1$
as a nonlocal vertex algebra is degenerate.
\el

\begin{proof} We have ${\bf 1}\in\ker \D$.
If $\dim \ker \D\ge 2$, there exists $a\in \ker \D$ such that
$a$ and ${\bf 1}$ are linearly independent. We have
$$0\ne a\otimes {\bf 1}\otimes 1-{\bf 1}\otimes a\otimes 1
\in V^{\otimes 2}\otimes \C((x_{1}))((x_{2}))$$
while
$$Z_{2}(a\otimes {\bf 1}\otimes 1-{\bf 1}\otimes a\otimes 1)=Y(a,x_{1}){\bf
1}-Y(a,x_{2}){\bf 1}=e^{x_{1}\D}a-e^{x_{2}\D}a=a-a=0.$$
Thus $Z_{2}$ is not injective.
\end{proof}

\bex{semidirect-sum}
{\em Here we give some examples of degenerate vertex algebras. 
Let $V$ be any vertex algebra and let $(W,Y_{W})$ be a $V$-module 
equipped with a linear operator $d$ on $W$ such that
$$[d,Y_{W}(v,x)]=\frac{d}{dx}Y_{W}(v,x)\ \ \ \mbox{ for }v\in V. $$
Extend the vertex operator map of $V$ to the vector space $V\oplus W$ by
$$Y(u+w,x)(v+w')=Y(u,x)v+Y_{W}(u,x)w'+e^{xd}Y_{W}(v,-x)w$$
for $u,v\in V,\; w,w'\in W$. 
Then (cf. \cite{li-form}, \cite{ll}) $V\oplus W$ is a vertex algebra 
with $V$ a vertex subalgebra and with $W$ an ideal.
Since $Y(w,x)w'=0$ for $w,w'\in W$,
if $W\ne 0$, $V\oplus W$ is degenerate.}
\eex

Let $V$ be a nonlocal vertex algebra. A {\em quasi $V$-module}
is a vector space $W$ equipped with a linear map
$Y_{W}$ from $V$ to $\E(W)$ such that
$Y_{W}({\bf 1},x)=1_{W}$ and such that for any $u,v\in V,\; w\in W$,
\begin{eqnarray}
p(x_{0}+x_{2},x_{2})Y_{W}(u,x_{0}+x_{2})Y_{W}(v,x_{2})w
=p(x_{0}+x_{2},x_{2})Y_{W}(Y(u,x_{0})v,x_{2})w
\end{eqnarray}
for some nonzero polynomial $p(x_{1},x_{2})$.

\bd{ddeg-g1va}
{\em Let $V$ be a nonlocal vertex algebra and let $(W,Y_{W})$ 
be be a quasi $V$-module. 
For each positive integer $n$, we define a linear map
$$\pi_{n}^{W}:V^{\otimes n}\otimes \C((x_{1}))\cdots ((x_{n}))\rightarrow 
\Hom (W,W((x_{1}))\cdots ((x_{n})))$$
  by
\begin{eqnarray}
\pi_{n}^{W}(u_{1}\otimes \cdots \otimes u_{n}\otimes f)=
f\cdot Y_{W}(u_{1},x_{1})\cdots Y_{W}(u_{n},x_{n}).
\end{eqnarray}}
\ed

\bl{ltwo-equivalent}
Let $V$ be a nonlocal vertex algebra and let $n$ be a positive integer.
If $Z_{n}$ is injective, then $\pi_{n}^{V}$ is injective.
If $V$ is a vertex algebra, the converse is also true.
\el

\begin{proof} The first assertion is clear.
Now, assume that $V$ is a vertex algebra. 
Then for any $u,v\in V$, there exists 
a nonnegative integer $k$ such that
$$(x_{1}-x_{2})^{k}Y(u,x_{1})Y(v,x_{2})
=(x_{1}-x_{2})^{k}Y(v,x_{2})Y(u,x_{1}).$$
In view of this, for any $a\in V$ and for any 
$A\in V^{\otimes n}\otimes \C((x_{1}))\cdots ((x_{n}))$, 
there exists a nonnegative integer $k$ such that
$$\left(\prod_{r=1}^{n}(x-x_{r})^{k}\right)Y(a,x)\pi_{n}^{V}(A)
=\left(\prod_{r=1}^{n}(x-x_{r})^{k}\right)\pi_{n}^{V}(A)Y(a,x).$$
If $A\in \ker Z_{n}$, we have
$$\left(\prod_{r=1}^{n}(x-x_{r})^{k}\right)\pi_{n}^{V}(A)Y(a,x){\bf 1}
=\left(\prod_{r=1}^{n}(x-x_{r})^{k}\right)Y(a,x)\pi_{n}^{V}(A){\bf 1}=0.$$
Thus $\pi_{n}^{V}(A)Y(a,x){\bf 1}=0$, 
which implies that $\pi_{n}^{V}(A)(a)=0$.
This proves that $A\in \ker \pi_{n}^{V}$. Therefore,
if $\pi_{n}^{V}$ is injective, so is $Z_{n}$.
\end{proof}

Let $\hat{\g}=\g\otimes \C[t,t^{-1}]\oplus \C c$ be the affine Lie algebra 
associated to a (not necessarily finite-dimensional) Lie algebra
$\g$ equipped with a nondegenerate symmetric invariant bilinear form
$\<\cdot,\cdot\>$, where $c$ is central and
\begin{eqnarray}
[a\otimes t^{m},b\otimes t^{n}]=[a,b]\otimes t^{m+n}+m\<a,b\>\delta_{m+n,0}c
\end{eqnarray}
for $a,b\in \g,\; m,n\in \Z$.
Associated to the
affine Lie algebra $\hat{\g}$ and any level $\ell\in \C$, we have a
vertex algebra $V_{\hat{\g}}(\ell,0)$, which is the induced
$\hat{\g}$-module from the $1$-dimensional $(\g[t]\oplus \C c)$-module
$\C_{\ell}$ $(=\C)$ with $\g[t]$ acting trivially and with $c$ acting
as the scalar $\ell$.  That is,
\begin{eqnarray}
V_{\hat{\g}}(\ell,0)=U(\hat{\g})\otimes _{U(\g[t]\oplus \C c)}\C_{\ell}
\simeq U(t^{-1}\g [t^{-1}])=S(t^{-1}\g[t^{-1}]).
\end{eqnarray}
This is an $\N$-graded
vertex algebra with $\C {\bf 1}$ being the degree zero subspace and
$\g$, identified with the degree $1$ subspace through the map
$a\mapsto a(-1){\bf 1}$ for $a\in \g$, generates $V_{\hat{\g}}(\ell,0)$
as a vertex algebra. If $\g$ is finite-dimensional and simple,
we often take $\<\cdot,\cdot\>$ to be the normalized Killing form
such that $\<\alpha,\alpha\>=2$ for long roots $\alpha$.
It is known (cf. \cite{fz}, \cite{dl}) that 
if $\ell$ is not equal to negative the dual Coxeter number of $\g$,
$V_{\hat{\g}}(\ell,0)$ is a vertex operator algebra (equipped with the
Segal-Sugawara conformal vector).
The following result was proved in \cite{ek}:

\bp{pek}
If $V_{\hat{\g}}(\ell,0)$ is an irreducible $\hat{\g}$-module, then
it is nondegenerate.
\ep

It was proved in \cite{li-simple} that every irreducible
vertex algebra of countable dimension over $\C$ is nondegenerate.
Next we generalize this result from vertex algebras to nonlocal vertex
algebras. First, we prove

\bl{lker-submodule}
Let $V$ be a nonlocal vertex algebra and 
let $(W,Y_{W})$ be any quasi $V$-module.
Let $n$ be a positive integer and consider 
$V^{\otimes n}\otimes \C((x_{1}))\cdots ((x_{n}))$
as a $V$-module with $V$ acting on the first factor $V$. Then
$\ker \pi_{n}^{W}$ is a $V$-submodule.
\el

\begin{proof} Suppose that
$$X=\sum_{i} u_{1i}\otimes \cdots 
\otimes u_{ni}\otimes f_{i}(x_{1},\dots,x_{n})
\in \ker \pi^{W}_{n}.$$
That is,
\begin{eqnarray}\label{eexp=0}
\sum_{i} f_{i}(x_{1},\dots,x_{n}) 
Y_{W}(u_{1i},x_{1})\cdots Y_{W}(u_{ni},x_{n})=0.
\end{eqnarray}
We must prove that for any $a\in V$ and any $w\in W$,
\begin{eqnarray}\label{eexp=1}
\sum_{i} f_{i}(x_{1},\dots,x_{n}) 
Y_{W}(Y(a,x_{0})u_{1i},x_{1})Y_{W}(u_{2i},x_{2})\cdots Y_{W}(u_{ni},x_{n})w=0.
\end{eqnarray}
Let $m_{2},\dots,m_{n}$ be arbitrarily fixed integers. 
Note that
$$\Res_{x_{2}}\cdots \Res_{x_{n}}x_{2}^{m_{2}}\cdots x_{n}^{m_{n}}
f_{i}(x_{1},\dots,x_{n}) Y_{W}(u_{2i},x_{2})\cdots Y_{W}(u_{ni},x_{n})w\in W\otimes \C((x_{1}))$$
for any $i$. Set
$$F_{i}=\Res_{x_{2}}\cdots \Res_{x_{n}}x_{2}^{m_{2}}\cdots x_{n}^{m_{n}}
f_{i}(x_{1},\dots,x_{n}) Y_{W}(u_{2i},x_{2})\cdots Y_{W}(u_{ni},x_{n})w.$$
For any $a\in V$, there exists $0\ne p(x_{1},x_{2})\in \C[x_{1},x_{2}]$ (depending on
$m_{2},\dots,m_{n}$) such that
\begin{eqnarray}
p(x_{0}+x_{1},x_{1})Y_{W}(a,x_{0}+x_{1})Y_{W}(u_{1i},x_{1})F_{i}
=p(x_{0}+x_{1},x_{1})Y_{W}(Y(a,x_{0})u_{1i}x_{1})F_{i}
\end{eqnarray}
for all the $i$. Then by applying
$p(x_{0}+x_{1},x_{1})Y_{W}(a,x_{0}+x_{1})
\Res_{x_{2}}\cdots \Res_{x_{n}}x_{2}^{m_{2}}\cdots x_{n}^{m_{n}}$
 to (\ref{eexp=0}) we get
\begin{eqnarray}
p(x_{0}+x_{1},x_{1})\sum_{i}Y_{W}(Y(a,x_{0})u_{1i},x_{1})F_{i}=0.
\end{eqnarray}
By cancellation we have
\begin{eqnarray}
\sum_{i}Y_{W}(Y(a,x_{0})u_{1i},x_{1})F_{i}=0.
\end{eqnarray}
That is, 
\begin{eqnarray}
& &\Res_{x_{2}}\cdots \Res_{x_{n}}x_{2}^{m_{2}}\cdots x_{n}^{m_{n}}
\sum_{i}f_{i}(x_{1},\dots,x_{n})Y_{W}(Y(a,x_{0})u_{1i},x_{1})Y_{W}(u_{2i},x_{2})\cdots
Y_{W}(u_{ni},x_{n})w\nonumber\\
& &=0.
\end{eqnarray}
Since $m_{2},\dots,m_{n}$ are arbitrary, we have (\ref{eexp=1}).
This proves that $\ker \pi^{W}_{n}$ is a $V$-submodule.
\end{proof}

\bt{tcountable}
Let $V$ be an irreducible nonlocal vertex algebra of 
countable dimension over $\C$ and let $(W,Y_{W})$ be any quasi $V$-module. 
Then for each positive integer $n$,
the linear map $\pi_{n}^{W}$ is injective.
\et

\begin{proof} Define $\pi_{0}^{W}$ to be the natural embedding of $\C$ into $\Hom (W,W)$.
We shall prove that if $\pi^{W}_{n}$ is not injective for some positive integer $n$,
then $\pi^{W}_{n-1}$ is not injective. Then the theorem follows immediately from induction. 

Assume that $\pi_{n}$ is not injective for some positive integer $n$.
By Lemma \ref{lker-submodule}, $\ker \pi^{W}_{n}$ 
is a $V$-submodule, where $V^{\otimes n}\otimes 
\C((x_{1}))\cdots ((x_{n}))$ is considered as a $V$-module with $V$ acting on the
first factor. 
Let $A$ be the associative algebra
generated by operators $v_{m}$ on $V^{\otimes n}\otimes 
\C((x_{1}))\cdots ((x_{n}))$ for $v\in V,\; m\in \Z$.
Since $V$ is of countable-dimensional, we have $\End_{A}V=\C$.
Let $0\ne X\in \ker \pi^{W}_{n}$.  We have
$$X=u_{1}\otimes e_{1}+\cdots +u_{r}\otimes e_{r},$$ 
where $u_{1},\dots,u_{r}\in V$ are linearly independent and
$0\ne e_{i}\in V^{\otimes (n-1)}\otimes \C((x_{1}))\cdots ((x_{n}))$.
By Jacobson density theorem
there exists $a\in A$ such that $au_{1}={\bf 1}$ and $au_{i}=0$ for
$1<i\le r$. In view of Lemma \ref{lker-submodule} we have
$$aX={\bf 1}\otimes e_{1}\in \ker \pi^{W}_{n}.$$
As ${\bf 1}\otimes e_{1}\in \ker \pi^{W}_{n}$, we have
$\Res_{x_{n}}x_{n}^{r}e_{1}\in \ker \pi^{W}_{n-1}$ for $r\in \Z$.
With $e_{1}\ne 0$, $\Res_{x_{n}}x_{n}^{r}e_{1}\ne 0$ for some $r\in \Z$.
Therefore, $\pi_{n-1}^{W}$ is not injective. This completes the proof.
\end{proof}

\bc{cnon-irreducible}
Every irreducible nonlocal vertex algebra of countable dimension 
over $\C$ is nondegenerate. In particular, simple vertex operator 
superalgebras over $\C$ (with finite-dimensional homogeneous
subspaces) are nondegenerate.
\ec

\br{rLaurent}
{\em Let $A$ be any unital (noncommutative) associative algebra 
and let $\partial$ be a derivation.
Then (\cite{bor}, cf. \cite{bk}, \cite{li-g1}) 
$A$ is a nonlocal vertex algebra with the identity $1$ as the vacuum vector 
and with
\begin{eqnarray}
Y(a,x)b=(e^{x\partial}a)b\;\;\;\mbox{ for }a,b\in A.
\end{eqnarray}
Let us denote this vertex algebra by $(A,\partial)$.
Now take $A=\C[t]$ and $\partial=d/dt$. We have
\begin{eqnarray}
Y(f(t),x)g(t)=\left(e^{x(d/dt)}f(t)\right) g(t)=f(t+x)g(t)
\;\;\;\mbox{ for }f(t),g(t)\in \C[t].
\end{eqnarray}
It is easy to see that the vertex algebra $(\C[t],d/dt)$ is simple, 
but the adjoint module is reducible. 
Since
\begin{eqnarray*}
& &Z_{2}\left(t\otimes 1\otimes (x_{1}-x_{2})^{-1}-1\otimes t\otimes 
(x_{1}-x_{2})^{-1}
-1\otimes 1\otimes 1\right)\nonumber\\
&=&(x_{1}-x_{2})^{-1}(t+x_{1})-(x_{1}-x_{2})^{-1}(t+x_{2})-1\\
&=&(x_{1}-x_{2})^{-1}(x_{1}-x_{2})-1\\
&=&0,
\end{eqnarray*}
$Z_{2}$ is not injective. This shows that
the vertex algebra $(\C[t],d/dt)$ is degenerate.
Thus the irreducibility assumption in Theorem \ref{tcountable}
is necessary and it cannot be replaced by the simplicity.
Furthermore, the degeneracy of vertex algebra $(\C[t],d/dt)$
implies that the vertex algebras $(\C[t,t^{-1}],d/dt)$
and $(\C((t)),d/dt)$ are degenerate as well. 
As $\C((t))$ is a field, $(\C((t)),d/dt)$ is irreducible.
Thus the countability assumption in Theorem \ref{tcountable}
is also necessary.
Also, note that even though it is degenerate, 
we have $Y(f(t),x)g(t)=f(t+x)g(t)\ne 0$ for any $0\ne f(t),g(t)\in \C((t))$.}
\er

\br{r}
{\em  For curiosity, one can show that $(\C[t^{-1}],d/dt)$ is nondegenerate.}
\er

Having proved that any irreducible nonlocal vertex algebra of countable dimension 
over $\C$ is nondegenerate, next we shall show that nonlocal vertex algebras
with a basis of P-B-W type
such as vertex algebras $V_{\hat{\g}}(\ell,0)$
are also nondegenerate. Our main tool will be 
increasing good filtrations which were studied in \cite{li-vpa}.

First we formulate the following straightforward result, 
which is classical in nature:

\bl{lgr-g1va}
Let $V$ be any nonlocal vertex algebra and let
$E=\{E_{n}\}_{n\in \Z}$ be an increasing filtration of $V$ 
satisfying the condition that ${\bf 1}\in E_{0}$ and 
\begin{eqnarray}\label{eincreasing-filtration-property}
u_{k}v\in E_{m+n}
\end{eqnarray}
for $u\in E_{m},\; v\in E_{n}$ and $k\in \Z$ with $m,n\in \Z$.
Form the $\Z$-graded vector space 
\begin{eqnarray}
\gr_{E}(V)=\coprod_{n\in \Z}E_{n}/E_{n-1}.
\end{eqnarray}
Then $\gr_{E}(V)$ is a nonlocal vertex algebra
with ${\bf 1}+E_{-1}$ as the vacuum vector and with
\begin{eqnarray}
(u+E_{m-1})_{k}(v+E_{n-1})=u_{k}v+E_{m+n-1}
\end{eqnarray}
for $u\in E_{m},\; v\in E_{n}$ with $m,n,k\in \Z$.
\el

We have:

\bp{pfiltration-non}
Let $V$ be a nonlocal vertex algebra. Assume that
there exists a lower truncated increasing filtration 
$E=\{E_{n}\}$ satisfying the property
(\ref{eincreasing-filtration-property}) such that the associated 
nonlocal vertex algebra $\gr_{E}(V)$ is nondegenerate. Then
$V$ is nondegenerate.
\ep

\begin{proof} First we construct a certain basis of $V$.
By assumption, there is an integer $r$ such that
$E_{n}=0$ for $n< r$. Let $\{ v^{\alpha}\;|\; \alpha\in I_{r}\}$ be a basis
of $E_{r}$, where $I_{r}$ is an index set. We extend this basis of $E_{r}$ to a basis
$\{ v^{\alpha}\;|\; \alpha\in I_{r}\}\cup
\{ v^{\alpha}\;|\; \alpha\in I_{r+1}\}$ of $E_{r+1}$, where $I_{r+1}$ is another index set,
and then extend this basis of $E_{r+1}$ to a basis of $E_{r+2}$, and so on.
In this way, we obtain a basis $\{ v^{\alpha}\;|\; \alpha \in I\}$
of $V$, where $I=\cup_{m\ge r}I_{m}$ (a disjoint union).
For $\alpha\in I_{m}$, we define $|\alpha|=m$. Then for $m\ge r$,
$$I_{m}=\{\alpha\in I\;|\; |\alpha|=m\}.$$
Let $n$ be a positive integer.
For ${\bf \alpha}=(\alpha_{1},\dots,\alpha_{n})\in I^{n}$, set
$$|\alpha|=|\alpha_{1}|+\cdots +|\alpha_{n}|.$$
{}From the construction of the basis $\{ v^{\alpha}\;|\; \alpha \in I\}$,
we have that $v^{\alpha}\in E_{|\alpha|}$ for $\alpha \in I$ and
that for any $m\in \Z$, $\{ v^{\alpha}+E_{m-1}\;|\; \alpha\in I_{m}\}$ is
a basis of $E_{m}/E_{m-1}$.

Assume that $Z_{n}$ for $V$ is not injective. Then
\begin{eqnarray}\label{eassumption-exp}
\sum_{{\bf \alpha} =(\alpha_{1},\dots,\alpha_{n})\in S} 
f_{(\alpha_{1},\dots,\alpha_{n})}(x_{1},\dots,x_{n}) 
Y(v^{\alpha_{1}},x_{1})\cdots Y(v^{\alpha_{n}},x_{n}){\bf 1}=0
\end{eqnarray}
for some finite subset $S$ of $I^{n}$ and 
$f_{\bf \alpha}(x_{1},\dots,x_{n})\ne 0$ for ${\bf \alpha}\in S$.
Set $k={\rm max}\{ |{\bf \alpha}|\;|\; {\bf \alpha}\in S\}$. 
Then $S=S_{r}\cup S_{r+1}\cup \cdots \cup S_{k}$.

{}From (\ref{eincreasing-filtration-property}),
for any ${\bf \alpha}=(\alpha_{1},\dots,\alpha_{n})\in S_{j}$,
we have
$$f_{(\alpha_{1},\dots,\alpha_{n})}(x_{1},\dots,x_{n}) 
Y(v^{\alpha_{1}},x_{1})\cdots Y(v^{\alpha_{n}},x_{n}){\bf 1}
\in E_{j}((x_{1}))\cdots ((x_{n})).$$
It follows from (\ref{eassumption-exp}) that
\begin{eqnarray*}
& &\sum_{{\bf \alpha}\in S_{k}}
f_{(\alpha_{1},\dots,\alpha_{n})}(x_{1},\dots,x_{n}) 
Y(v^{\alpha_{1}},x_{1})\cdots Y(v^{\alpha_{n}},x_{n}){\bf 1}\\
&=&-\sum_{m=r}^{k-1}\sum_{{\bf \alpha}\in S_{m}}
f_{(\alpha_{1},\dots,\alpha_{n})}(x_{1},\dots,x_{n}) 
Y(v^{\alpha_{1}},x_{1})\cdots Y(v^{\alpha_{n}},x_{n}){\bf 1}\\
&\in& E_{k-1}((x_{1}))\cdots ((x_{n})).
\end{eqnarray*}
Thus
$$\sum_{{\bf \alpha}\in S_{k}}
f_{(\alpha_{1},\dots,\alpha_{n})}(x_{1},\dots,x_{n}) 
Y(v^{\alpha_{1}}+E_{|\alpha_{1}|-1},x_{1})\cdots 
Y(v^{\alpha_{n}}+E_{|\alpha_{n}|-1},x_{n}){\bf 1}=0$$
in $(\gr_{E}(V))^{\otimes n}\otimes \C((x_{1}))\cdots ((x_{n}))$.
{}From the construction of $\{ v^{\alpha}\;|\; \alpha\in I\}$,  we have
$$\sum_{{\bf \alpha}\in S_{k}}
(v^{\alpha_{1}}+E_{|\alpha_{1}|-1})\otimes \cdots \otimes 
(v^{\alpha_{n}}+E_{|\alpha_{n}|-1})\otimes 
f_{(\alpha_{1},\dots,\alpha_{n})}(x_{1},\dots,x_{n}) 
\ne 0$$
in $(\gr_{E}(V))^{\otimes n}\otimes \C((x_{1}))\cdots ((x_{n}))$.
This contradicts with our assumption that
$\gr_{E}(V)$ is nondegenerate.
\end{proof}

Just as with classical associative algebras, 
for nonlocal vertex algebras we can obtain
increasing filtrations from any generating subsets.

\bp{pgenerating-filtration}
Let $V$ be a nonlocal vertex algebra and let $U$ be 
any generating subset of $V$.
For any nonnegative integer $n$, define 
\begin{eqnarray}
E_{n}^{U}={\rm span}\{ u^{(1)}_{m_{1}}\cdots u^{(r)}_{m_{r}}{\bf 1}
\;|\; 0\le r\le n,\; u^{(i)}\in U,m_{i}\in \Z\}. 
\end{eqnarray}
Then $\{E_{n}^{U}\}_{n\in \N}$ is an increasing filtration of $V$ 
satisfying the condition (\ref{eincreasing-filtration-property}).
Furthermore, if for $u,v\in U$, $m,n\in \Z$, and $k\ge 0$,
\begin{eqnarray}
(u_{m}v_{n}-v_{n}u_{m})E_{k}\subset E_{k+1},
\end{eqnarray}
then the associated nonlocal vertex algebra $\gr_{E}(V)$ is 
a (strictly) commutative vertex algebra.
\ep

\begin{proof} As $U$ generates $V$ as a nonlocal vertex algebra, 
$\{E_{n}^{U}\}_{n\in \N}$ is an increasing filtration of $V$. 
To establish the property (\ref{eincreasing-filtration-property}),
we are going to use induction to prove that for any $m\ge 0$,
$$a_{k}E_{n}\subset E_{m+n}$$
for all $a\in E_{m},\; k\in \Z$. Since $E_{0}=\C {\bf 1}$,
it is true for $m=0$. Now assume $a\in E_{m+1}$.
Notice that for any $m\ge 0$,
$$E_{m+1}={\rm span}\{ u_{r}v\;|\; u\in U,\; r\in \Z,\; v\in E_{m}\}.$$
In view of this, it suffices to consider $a=u_{r}v$ 
for $u\in U,\; r\in \Z,\; v\in E_{m}$.
For $w\in E_{n}$, there exists a nonnegative integer $l$ such that
$$(x_{0}+x_{2})^{l}Y(Y(u,x_{0})v,x_{2})w
=(x_{0}+x_{2})^{l}Y(u,x_{0}+x_{2})Y(v,x_{2})w.$$
With $v\in E_{m},\; w\in E_{n}$, using inductive hypothesis we have
$$Y(u,x_{0}+x_{2})Y(v,x_{2})w\in Y(u,x_{0}+x_{2})E_{n+m}((x_{2}))
\subset E_{n+m+1}((x_{0}))((x_{2})).$$
Consequently we have
$$(x_{0}+x_{2})^{l}Y(Y(u,x_{0})v,x_{2})w
\in E_{n+m+1}((x_{0}))((x_{2})).$$
Since $Y(Y(u,x_{0})v,x_{2})w\in V((x_{2}))((x_{0}))$, we have
$$Y(Y(u,x_{0})v,x_{2})w
\in E_{n+m+1}((x_{2}))((x_{0})).$$
In particular, we have $a_{k}E_{n}\subset E_{n+m+1}$.

Notice that $U\subset E_{1}$. As $U$ generates $V$ as a nonlocal vertex algebra,
$U/E_{0}$ generates $\gr_{E}(V)$ as a nonlocal vertex algebra. Then the last assertion is clear.
\end{proof}

We next recall certain results on constructing increasing filtrations
for vertex algebras from \cite{li-vpa}:

\bl{lvpa}
Let $V$ be an (ordinary) vertex algebra. 
Let $E=\{ E_{n}\}$ be an increasing good
filtration of $V$ in the sense that 
${\bf 1}\in E_{0}$ and that for $u\in E_{m},\; v\in E_{n}$
$$u_{k}v\in E_{m+n}\;\;\;\mbox{ for any }k\in \Z,$$ 
and $u_{k}v\in E_{m+n-1}$ for $k\ge 0$. 
Then $\gr_{E}(V)$ is a strictly commutative vertex algebra, 
which amounts to a unital commutative associative differential algebra.
\el

\bt{tgeneral-pbw-spanning-property}
Let $V$ be an (ordinary) vertex algebra and
let $U$ be a subspace of $V$ equipped with a vector space decomposition
$U=\coprod_{n\ge 1}U_{n}$. For $n\ge 0$, denote by $E_{n}^{U}$ the subspace of $V$
linearly spanned by the vectors
$$u^{(1)}_{-k_{1}}\cdots u^{(r)}_{-k_{r}}{\bf 1}$$
for $r\ge 0,\; u^{(i)}\in U,\; k_{i}\ge 1$
with $\wt u^{(1)}+\cdots +\wt u^{(r)}\le n$. 
Assume that $U$ generates $V$ as a vertex algebra and that
for $u\in U_{m},\; v\in U_{n}$ with $m,n\ge 1$,
$u_{i}v\in E_{m+n-1}$ for all $i\ge 0$.  
Then $\{E^{U}_{n}\}$ is an increasing good filtration of $V$ and
$\gr_{E^{U}}(V)$ is a (strictly) commutative vertex algebra.
Furthermore, any good filtration $\{E_{n}\}$ of $V$ with $E_{n}=0$ for $n<0$ 
and with $E_{0}=\C {\bf 1}$ can be obtained this way.
\et

Next, we establish nondegeneracy
for better known vertex operator algebras.
Let $U$ be any vector space over $\C$.
Set
\begin{eqnarray}
FCD(U)=S(\C[\partial]\otimes U),
\end{eqnarray}
the symmetric algebra over the vector space $\C[\partial]\otimes U$,
where $\partial$ is a formal variable.  Then $FCD(U)$ is a commutative
associative algebra equipped with a derivation $\partial$.  We call
$FCD(U)$ the {\em free commutative differential algebra over $U$}.

\bp{pcomm-diff}
Let $U$ be any vector space of countable dimension over $\C$. 
The free commutative differential algebra $FCD(U)$ on $U$,
viewed as a vertex algebra, is nondegenerate.
\ep

\begin{proof} Pick up a basis $\{ u_{\alpha}\;|\; \alpha\in I\}$ of $U$.
Let $U^{'}$ be the subspace of $U^{*}$, linearly spanned by
the linear functionals $u_{\alpha}^{*}$ dual to $u_{\alpha}$ for $\alpha\in I$.
Set $\h=U\oplus U^{'}$, a vector space of countable dimension over $\C$.
Equip $\h$ with the standard (nondegenerate) symmetric bilinear form $\<\cdot,\cdot\>$
$$\<a+a',b+b'\>=\<a,b'\>+\<a',b\>$$ for $a,b\in U,\; a',b'\in
U^{'}$. 
As $V_{\hat{\h}}(1,0)$ is an irreducible vertex algebra 
of countable dimension over $\C$,
by Corollary \ref{cnon-irreducible}, $V_{\hat{\h}}(1,0)$ is
nondegenerate.  Consider the vertex subalgebra $\<U\>$ of
$V_{\hat{\h}}(1,0)$, generated by the subspace $U$.  
In view of Lemma \ref{lnon-subva}, $\<U\>$ is nondegenerate.
For $a,b\in U$ and for $m,n\in \Z$, we have
$$[a(m),b(n)]\;(=[a,b](m+n)+m\<a,b\>\delta_{m+n,0})=0.$$ 
As $U$ generates $\<U\>$ as a vertex algebra,
$\<U\>$ is a (strictly) commutative vertex algebra. Equivalently,
$\<U\>$ is a commutative differential algebra. In view of the P-B-W theorem,
we have $\<U\>=S(t^{-1}U[t^{-1}])$ as a vector space.
Consequently, we have that $FCD(U)=\<U\>$ 
as a commutative differential algebra and then
as a vertex algebra.  Therefore,
$FCD(U)$ is nondegenerate.
\end{proof}

Combining Propositions \ref{pfiltration-non} and \ref{pcomm-diff} 
we immediately have:

\bt{tpbw-va-verma}
Let $V$ be a nonlocal vertex algebra of countable dimension over $\C$. 
Assume that there exists an increasing filtration $E=\{E_{n}\}_{n\ge 0}$ of $V$, 
satisfying the condition (\ref{eincreasing-filtration-property}),
 such that the associated nonlocal vertex algebra
$\gr_{E}(V)$ is isomorphic to the free commutative differential algebra
$FCD(U)$ over a vector space $U$. 
Then $V$ is nondegenerate.
\et

In view of Theorem \ref{tweak-strongqva} we immediately have:

\bc{cqva}
Let $V$ be a weak quantum vertex algebra of countable dimension over $\C$.
Assume that there exists an increasing filtration 
$E=\{E_{n}\}_{n\ge 0}$ of $V$, 
satisfying the condition (\ref{eincreasing-filtration-property}),
such that the associated nonlocal vertex algebra
$\gr_{E}(V)$ is isomorphic to the free commutative differential algebra
$FCD(U)$ over a vector space $U$. Then $V$ is 
a nondegenerate quantum vertex algebra.
\ec

Now we are ready to establish the nondegeneracy of vertex algebra $V_{\hat{\g}}(\ell,0)$
in the general case (cf. \cite{ek}).

\bp{paffine-verma}
Let $\g$ be any Lie algebra of countable dimension over $\C$ equipped with 
a nondegenerate symmetric invariant bilinear form.
Then for any complex number $\ell$, the vertex algebra
$V_{\hat{\g}}(\ell,0)$ is nondegenerate.
\ep

\begin{proof} Recall that $V_{\hat{\g}}(\ell,0)$ as a vertex algebra 
is generated by $\g$. Taking $U=U_{1}=\g$ in Theorem \ref{tgeneral-pbw-spanning-property}
we obtain a good increasing filtration $E^{\g}=\{E_{n}^{\g}\}_{n\ge 0}$ 
of $V_{\hat{\g}}(\ell,0)$.
For any positive integer $n$,
$E^{\g}_{n}$ is the linear span of the vectors 
$$a^{(1)}(-m_{1})\cdots a^{(r)}(-m_{r}){\bf 1}$$
for $0\le r\le n,\; a^{(1)},\dots,a^{(r)}\in \g,\; m_{1},\dots,m_{r}\ge 1$.
It follows from the P-B-W theorem that $V_{\hat{\g}}(\ell,0)=U(t^{-1}\g[t^{-1}])$ 
as a vector space. The filtration $E$ coincides with the standard filtration
on $U(t^{-1}\g[t^{-1}])$. Then $\gr_{E}(V_{\hat{\g}}(\ell,0))=S(t^{-1}\g[t^{-1}])$.
Thus $\gr_{E}(V_{\hat{\g}}(\ell,0)))$, a (strictly)  commutative vertex algebra,  
is isomorphic to the free commutative differential algebra
over $\g$. By Theorem \ref{tpbw-va-verma}, $V_{\hat{\g}}(\ell,0)$ is nondegenerate.
\end{proof}

Consider the Virasoro Lie algebra $Vir=\sum_{n\in\Z}\C L(n)\oplus \C c$.
Recall (cf. \cite{fz}) that for any complex number $\ell$, we have a
vertex operator algebra $V_{Vir}(\ell,0)$, associated to the Virasoro Lie algebra $Vir$
of central charge $\ell$. By definition, 
\begin{eqnarray}
V_{Vir}(\ell,0)=U(Vir)\otimes_{U(Vir_{(\ge 0)})}\C_{\ell},
\end{eqnarray}
where $\C_{\ell}$ is the $1$-dimensional module for the Lie subalgebra
$$Vir_{(\ge 0)}=\sum_{n\ge 0}\C L(n-1)\oplus \C c$$
 with $c$ acting as scalar $\ell$
and with $\sum_{n\ge 0}\C L(n-1)$ acting as zero.

\bp{pviraso-verma}
For any complex number $\ell$, the vertex operator algebra
$V_{Vir}(\ell,0)$ is nondegenerate.
\ep

Let ${\mathcal{L}}$ be a Lie algebra equipped with a derivation $d$.
Then $d$ naturally extends to a derivation of $U({\mathcal{L}})$, 
the universal enveloping algebra of ${\mathcal{L}}$. 
{}From Remark \ref{rLaurent}, $U({\mathcal{L}})$ is naturally a nonlocal vertex algebra. 

\bp{paffine-like-half}
Let ${\mathcal{L}}$ be a Lie algebra equipped with a derivation
$d$ such that ${\mathcal{L}}$ is a free $\C[d]$-module with a countable base.
Then the nonlocal vertex algebra
$U({\mathcal{L}})$ is nondegenerate.
\ep

\begin{proof} Recall that for $a,b\in U({\mathcal{L}})$ we have $Y(a,x)b=(e^{xd}a)b$.
That is, $a_{n}b=0$ for $n\ge 0$ and $a_{n}b=\frac{1}{(-n-1)!}(d^{-n-1}a)b$ for $n<0$.
As ${\mathcal{L}}$ generates $U({\mathcal{L}})$ as an algebra,
${\mathcal{L}}$ generates $U({\mathcal{L}})$ as a nonlocal vertex algebra.
Also, for the increasing filtration $E$ defined in Proposition \ref{pgenerating-filtration} we have 
$$E_{n}={\rm span}\{ a^{(1)}\cdots a^{(r)}\;|\; 0\le r\le n,\; a^{(i)}\in {\mathcal{L}}\}$$
for $n\ge 0$. Thus, the filtration $E$ is exactly the standard filtration for the universal enveloping algebra
$U({\mathcal{L}})$. It was well known that $\gr_{E}(U({\mathcal{L}}))=S({\mathcal{L}})$
(the symmetric algebra).
By Theorem \ref{tpbw-va-verma} the nonlocal vertex algebra $U({\mathcal{L}})$ is 
nondegenerate.
\end{proof}

We next consider current Lie algebras.
Let $\g$ be any Lie algebra. With $t^{-1}\C[t^{-1}]$ being a commutative associative algebra,
$\g \otimes t^{-1}\C [t^{-1}]$ is naturally a Lie algebra. Furthermore, $L(-1)=1\otimes d/dt$ is
a derivation. It is clear that $\g \otimes t^{-1}\C [t^{-1}]$ is a free $\C[L(-1)]$-module.
In view of Proposition \ref{paffine-like-half} we immediately have:

\bc{caffine-half}
Let $\g$ be any Lie algebra of countable dimension over $\C$. 
Then the nonlocal vertex algebra $U(\g\otimes t^{-1}\C[t^{-1}])$
associated to the noncommutative differential algebra structure
is nondegenerate.
\ec

Furthermore we have:

\bp{paffine-half-quantumva}
Let $\g$ be any Lie algebra of countable dimension over $\C$. Then the nonlocal vertex algebra
$U(\g\otimes t^{-1}\C[t^{-1}])$ associated to the noncommutative differential algebra
structure is a nondegenerate quantum vertex algebra. 
\ep

\begin{proof} For $a\in \g$, set
\begin{eqnarray}
a(x)^{+}=\sum_{n\ge 1}(a\otimes t^{-n})x^{n-1}\in (\g\otimes t^{-1}\C[t^{-1}])[[x]].
\end{eqnarray}
We have
$$Y(a\otimes t^{-1},x)=e^{xL(-1)}(a\otimes t^{-1})=a\otimes (t+x)^{-1}
=\sum_{n\ge 0}(a\otimes t^{-n-1})x^{n}=a(x)^{+}.$$
It is clear that $\g\otimes t^{-1}$ generates $U(\g\otimes t^{-1}\C[t^{-1}])$ 
as a nonlocal vertex algebra.
For $a,b\in \g$, by a straightforward calculation we have
\begin{eqnarray}\label{e-comm-a+b+}
[a(x_{1})^{+},b(x_{2})^{+}]=(x_{2}-x_{1})^{-1}\left([a,b](x_{2})^{+}-[a,b](x_{1})^{+}\right).
\end{eqnarray}
{}From this, we have
$$(x_{1}-x_{2})a(x_{1})^{+}b(x_{2})^{+}
=(x_{1}-x_{2})b(x_{2})^{+}a(x_{1})^{+}+[a,b](x_{2})^{+}\cdot 1-1\cdot [a,b](x_{1})^{+}.
$$
This shows that $\{ Y(a\otimes t^{-1},x)=a(x)^{+}\;|\; a\in \g\}$
is $\S$-local and that $\g\otimes t^{-1}$ is an $\S$-local generating subspace
of $U(\g\otimes t^{-1}\C[t^{-1}])$.
By Lemma \ref{lS-local},
$U(\g\otimes t^{-1}\C[t^{-1}])$ is a weak quantum vertex algebra.
As it has been proved in Corollary \ref{caffine-half} that $U(\g\otimes t^{-1}\C[t^{-1}])$ 
is nondegenerate, the rest follows from Theorem \ref{tweak-strongqva}.
\end{proof}

\br{rqyb-operator-appl}
{\em Let $\g$ be a Lie algebra of countable dimension over $\C$.
Set $U=\g\oplus \C$.
It follows from Proposition \ref{paffine-half-quantumva}, 
we have a unitary rational quantum Yang-Baxter operator $\S$ on $U$ such that
\begin{eqnarray}
& &\S(x)(1\otimes 1)=1\otimes 1,\\
& &\S(x)(1\otimes a)=1\otimes a,\\
& &\S(x)(a\otimes 1)=a\otimes 1,\\
& &\S(x)(a\otimes b)=a\otimes b+x^{-1}({\bf 1}\otimes [a,b]-[a,b]\otimes {\bf 1})
\end{eqnarray}
for $a,b\in \g$. }
\er

\br{rsingular-part}
{\em  Consider the Lie algebra $\g[t]$.
Let $d$ be any derivation of the polynomial algebra $\C[t]$.
Then $d$ is naturally a derivation of the Lie algebra
$\g[t]$, hence a derivation of $U(\g[t])$.
Now we have a nonlocal vertex algebra $U(\g[t])$.
But, since $\D=d$ and $\g\oplus \C\subset \ker d$, 
by Lemma \ref{lnon-kernel} $U(\g[t])$ is degenerate.}
\er

\br{ryangian}
{\em Note that Yangian $Y(\g)$ is closely related to the Lie algebra
$\g[t]$. In fact, the associated graded algebra of
$Y(\g)$ with respect to a certain increasing filtration 
is isomorphic to $U(\g[t])$. 
We shall study Yangian $Y(\g)$ in terms of quantum vertex algebras 
in a separate publication.}
\er

Next we construct some examples of nondegenerate quantum vertex algebras
{}from certain Lie algebras. Consider a semidirect product Lie algebra $K\oplus \g$, where
$K$ is an ideal and $\g$ is a Lie subalgebra.
Suppose that $\<\cdot,\cdot\>$ is a symmetric invariant bilinear form on $K$.
Extend $\<\cdot,\cdot\>$ to a bilinear form on $K\oplus \g$ by
$\< u+a,v+b\>=\<u,v\>$ for $u,v\in K,\; a,b\in \g$.
We have an affine Lie algebra associated to $(K\oplus \g,\<\cdot,\cdot\>)$.
Consider the following Lie subalgebra
\begin{eqnarray}
{\mathcal{L}}=K\otimes \C[t,t^{-1}]\oplus \C c\oplus (\g\otimes t^{-1}\C[t^{-1}]),
\end{eqnarray}
which is the semidirect product Lie algebra of $\hat{K}$ 
with $\g\otimes t^{-1}\C[t^{-1}]$.
Let $\ell$ be any complex number. Let $\C_{\ell}$ be the $1$-dimensional module
for $K[t]\oplus \C c$ with $K[t]$ acting as zero and with $c$ acting as scalar $\ell$.
Form the induced $\mathcal{L}$-module
\begin{eqnarray}
V_{\mathcal{L}}(\ell,0)=U({\mathcal{L}})\otimes_{U(K[t]\oplus \C c)} \C_{\ell}.
\end{eqnarray}
In view of the P-B-W theorem, we have
\begin{eqnarray}
V_{\mathcal{L}}(\ell,0)=U((K\oplus \g)\otimes t^{-1}\C[t^{-1}]),
\end{eqnarray}
as a vector space. Embed $K\oplus \g$ into $V_{\mathcal{L}}(\ell,0)$
through the map $u+a\mapsto u(-1){\bf 1}+a(-1){\bf 1}$.

\bp{psemiproduct}
There exists a unique weak quantum vertex algebra structure
on $V_{\mathcal{L}}(\ell,0)$ such that 
\begin{eqnarray}
Y(u+a,x)=u(x)+a(x)^{+}\;\;\;\mbox{ for }u\in K,\; a\in \g.
\end{eqnarray}
Furthermore, $V_{\mathcal{L}}(\ell,0)$ is a nondegenerate
quantum vertex algebra.
\ep

\begin{proof} We shall apply Theorem \ref{tweak-qva-construction}.
Set
$$U=\{ u(x),a(x)^{+}\;|\; u\in K,\; a\in \g\}.$$ We now prove that 
$U$ is $\S$-local.
For $u,v\in K$, we have
$$[u(x_{1}),v(x_{2})]=[u,v](x_{2})x_{1}^{-1}\delta\left(\frac{x_{2}}{x_{1}}\right)
+\ell \<u,v\> \frac{\partial}{\partial x_{2}}x_{1}^{-1}\delta\left(\frac{x_{2}}{x_{1}}\right).$$
For $u\in K,\; a\in \g$, we have
$$[a(x_{1})^{+},u(x_{2})]=[a,u](x_{2})(x_{2}-x_{1})^{-1}.$$
For $a,b\in \g$, we have
$$[a(x_{1})^{+},b(x_{2})^{+}]=
(x_{2}-x_{1})^{-1}([a,b](x_{2})^{+}-[a,b](x_{1})^{+}).$$
It follows that $U$ is $\S$-local. By Theorem \ref{tX}, $U$ generates 
a weak quantum vertex algebra $\<U\>$ inside $\E(V)$. It follows from
Proposition \ref{pY} that $\<U\>$ is naturally 
an ${\mathcal{L}}$-module of level $\ell$. Then there exists 
a unique ${\mathcal{L}}$-homomorphism $\psi$ from $V$ to $\<U\>$, sending 
$1_{W}$ to ${\bf 1}$. Now, we have all the conditions 
assumed in Theorem \ref{tweak-qva-construction}.
\end{proof}

Notice that vertex algebras are analogous to unital commutative associative algebras.
It is a well known fact that any unital commutative associative algebra
without zero divisor can be embedded into a field, namely a simple (and irreducible) 
commutative associative algebra. In view of this we formulate the following:

\bconj{embedding}
Every  nondegenerate vertex algebra can be embedded into a 
simple (irreducible) vertex algebra. 
\econj

\section{Quantum vertex algebras associated with 
associative algebras of Zamolodchikov-Faddeev type}
In this section we construct a family of quantum vertex algebras 
closely related to Zamolodchikov-Faddeev algebras.

Let $H$ be a vector space equipped with a bilinear form $\<\cdot,\cdot\>$ 
and let $\S(x): H\otimes H\rightarrow H\otimes H\otimes \C[[x]]$ be a 
linear map. Recall from the end of Section 2 
the notions of $(H,\S)$-module,
vacuum-like vector and vacuum $(H,\S)$-module.

Set $\S(x)=\sum_{n\ge 0}S_{n}x^{n}$ where $S_{n}\in \End (H\otimes H)$.
Denote by $A(H,S_{0})$ the quotient algebra of
the (free) tensor algebra $T(H\otimes t^{-1}\C[t^{-1}])$ modulo the following relations:
\begin{eqnarray}
(a\otimes t^{-m})(b\otimes t^{-n})-\sum_{i=1}^{r} (b^{(r)}\otimes t^{-n})(a^{(r)}\otimes t^{-m})
=0
\end{eqnarray}
for $a,b\in H,\; m,n\ge 1$, where $S_{0}(b\otimes a)=\sum_{i=1}^{r}b^{(r)}\otimes a^{(r)}$.
Notice that $A(H,1)$ is the symmetric algebra over $H\otimes t^{-1}\C[t^{-1}]$
while $A(H,-1)$ is the exterior algebra.

\bl{lpbw-span-property}
Let $H$ be a vector space equipped with a bilinear form
$\<\cdot,\cdot\>$ and let $S(x)$ be a linear map from
$H\otimes H$ to $H\otimes H\otimes \C[[x]]$.
Let $W$ be any vacuum $(H,\S)$-module with a vacuum-like vector $e$ as its generator.
For $k\ge 0$, let $W[k]$ be the span of the vectors
\begin{eqnarray}
a^{(1)}(-m_{1})\cdots a^{(r)}(-m_{r})e
\end{eqnarray}
for $r\ge 0,\; a^{(1)},\dots,a^{(r)}\in H,\;m_{i}\ge 1$ with $m_{1}+\cdots
+m_{r}\le k$.
Then $W[k]$ for $k\ge 0$ form an increasing filtration of $W$ such that
for $a\in H,\; m,k\in \Z$,
\begin{eqnarray}
a(m)W[k]&\subset& W[k-m],\label{eamwk-gen}\\
a(m)W[k]&\subset& W[k-m-1]\;\;\;\mbox{ whenever }m\ge 0,\label{eamwk-nonneg}
\end{eqnarray}
where by convention we set $W[k]=0$ for $k<0$.
Furthermore, $\gr W=\coprod_{k\ge 0}W[k]/W[k-1]$ is an $A(H,S_{0})$-module with
$e\in W[0]=W[0]/W[-1]$ as a generator, where
 $$a(-m)(w+W[k-1])=a(-m)w+W[k+m-1]\in W[k+m]/W[k+m-1]$$
for $a\in H,\; m\ge 1,\; w+W[k-1]\in W[k]/W[k-1]$ with $k\ge 0$.
\el

\begin{proof} {}From definition, we have $W[0]=\C e$ and
by convention $W[k]=0$ for $k<0$.
By definition we have
\begin{eqnarray}\label{efiltration-property-negatice-new}
a(m)W[k]\subset W[k-m] \;\;\;\mbox{ for }a\in H,\; m<0.
\end{eqnarray}
We next prove that
\begin{eqnarray}\label{efiltration-property-new}
a(m)W[k]\subset W[k-m-1] \;\;\;\mbox{ for }a\in H,\; m\ge 0.
\end{eqnarray}

For $a,b\in H$, as $\S(x)(b\otimes a)\in H\otimes H\otimes \C[[x]]$, 
we have
$$\iota_{x_{2},x_{1}}\S(x_{2}-x_{1})(b\otimes a)
=\sum_{r=1}^{s}\sum_{i\ge 0}\sum_{j=0}^{i}
\beta_{rij}(b^{(r)}\otimes a^{(r)}\otimes x_{2}^{i-j}x_{1}^{j}),$$
where $\beta_{rij}\in \C,\; b^{(i)},a^{(i)}\in H$. 
For $m,n\in \Z,\;w\in W$, {}from relation (\ref{emain-relation-qLie}) we have
\begin{eqnarray}
a(m)b(n)w=\sum_{r=1}^{s}\sum_{i\ge 0}\sum_{j=0}^{i}
\beta_{rij}b^{(r)}(n+i-j)a^{(r)}(m+j)w
+\<a,b\>\delta_{m+n+1,0}w.
\label{eexplicituv-new}
\end{eqnarray}
As $W[0]=\C e$, from the assumption we have
$$a(m)W[0]=0=W[-m-1]\;\;\;\mbox{ for }a\in H,\; m\ge 0.$$
Now, (\ref{efiltration-property-new}) follows immediately from induction (on $k$) and 
(\ref{eexplicituv-new}).

{}From (\ref{efiltration-property-negatice-new}) 
and (\ref{efiltration-property-new}),
$\cup_{k\ge 0}W[k]$ is
closed under the action of $a(m)$ for $a\in H,\; m\in \Z$.
As $W=T(H\otimes \C[t,t^{-1}])e$, we have $W=\cup_{k\ge 0}W[k]$. 
Thus $W[k]$ for $k\ge 0$
form an increasing filtration of $W$.
Using (\ref{eexplicituv-new}), 
(\ref{efiltration-property-negatice-new}) and
(\ref{efiltration-property-new}) we have
\begin{eqnarray*}
& &\left(a(-m)b(-n)-\sum_{r=1}^{s}
\beta_{r00}b^{(r)}(-n)a^{(r)}(-m)\right)W[k]\\
&=&\sum_{r=1}^{s}\sum_{i\ge 1}\sum_{j=0}^{i}
\beta_{rij}b^{(r)}(-n+i-j)a^{(r)}(-m+j)W[k]\\
&\subset& W[k+m+n-1]
\end{eqnarray*}
for $m,n\ge 1$. Notice that
$$S_{0}(b\otimes a)=\sum_{r=1}^{s}\beta_{r00}(b^{(r)}\otimes a^{(r)}).$$
Now, the second assertion is clear.
\end{proof}

We also have:

\bc{cpbw-span-property}
Let $H,\<\cdot,\cdot\>,\S$ and $W,e$ be given as in Lemma \ref{lpbw-span-property}.
Then
\begin{eqnarray}
a^{(1)}(n_{1})\cdots a^{(r)}(n_{r})e=0
\label{etruncation-vk-new}
\end{eqnarray}
whenever $n_{1}+\cdots +n_{r}\ge 0$ for $r\ge 1,\; a^{(i)}\in H,\;n_{i}\in \Z$.
\ec

\begin{proof} If $n_{1}+\cdots +n_{r}\ge 0$, 
we have $n_{i}\ge 0$ for some $i$.
Using (\ref{eamwk-gen}) and (\ref{eamwk-nonneg}) 
we have 
$$a^{(1)}(n_{1})\cdots a^{(r)}(n_{r})e\in W[-n_{1}-\cdots -n_{r}-1]=0,$$
proving (\ref{etruncation-vk-new}).
\end{proof}

Let $H,\<\cdot,\cdot\>$, and $\S$ be given as in Lemma \ref{lpbw-span-property}.
Let $T(H\otimes \C[t,t^{-1}])^{+}$
be the subspace spanned by the vectors
$$(a^{(1)}\otimes t^{n_{1}})\cdots (a^{(r)}\otimes t^{n_{r}})$$
for $r\ge 1,\; a^{(i)}\in H,\;n_{i}\in \Z$ with
$n_{1}+\cdots +n_{r}\ge 0$.
Set
\begin{eqnarray}
J=T(H\otimes \C[t,t^{-1}])T(H\otimes \C[t,t^{-1}])^{+},
\end{eqnarray}
a left ideal of $T(H\otimes \C[t,t^{-1}])$. We then set
\begin{eqnarray}
\tilde{V}(H,\S)=T(H\otimes \C[t,t^{-1}])/J,
\end{eqnarray}
a left $T(H\otimes \C[t,t^{-1}])$-module. Clearly, $\tilde{V}(H,\S)$ is cyclic on the vector
$\tilde{1}=1+J$ where $a(n)\tilde{1}=0$ for $a\in H,\; n\ge 0$. 
Furthermore, for $a\in H,\; w\in \tilde{V}(H,\S)$, we have
\begin{eqnarray}
a(m)w=0\;\;\;\mbox{ for $m$ sufficiently large}.
\end{eqnarray}
Now we define $V(H,\S)$ to be the quotient $T(H\otimes \C[t,t^{-1}])$-module
of $\tilde{V}(H,\S)$ modulo the following relations:
\begin{eqnarray}
a(x_{1})b(x_{2})w
-\sum_{i=1}^{r}\iota_{x_{2},x_{1}}(f_{i}(x_{2}-x_{1}))b^{(i)}(x_{2})a^{(i)}(x_{1})w
=x_{2}^{-1}\delta\left(\frac{x_{1}}{x_{2}}\right)\<a,b\>w
\end{eqnarray}
for $a,b\in H,\; w\in \tilde{V}(H,\S)$, where 
$\S(x)(b\otimes a)=\sum_{i=1}^{r}b^{(i)}\otimes a^{(i)}\otimes f_{i}(x)$. 
Denote by ${\bf 1}$ the image of $\tilde{1}$ in $V(H,\S)$. From definition,
$V(H,\S)$ is a vacuum $(H,\S)$-module with vacuum-like vector ${\bf 1}$ as its generator.
In view of Corollary \ref{cpbw-span-property} we immediately have:

\bp{puniversal}
Let $H$ be a vector space equipped with a bilinear form $\<\cdot,\cdot\>$
and let $\S(x): H\otimes H\rightarrow H\otimes H\otimes \C[[x]]$
be a linear map. The above constructed $(H,\S)$-module $V(H,\S)$ is 
a universal vacuum $(H,\S)$-module.
\ep	

Combining Propositions \ref{puniversal} and \ref{pqLie-wqva} and 
Corollary \ref{cqva} we immediately obtain:

\bt{tzd-qva}
Let $H$ be a finite-dimensional vector space over $\C$ 
equipped with a bilinear form
$\<\cdot,\cdot\>$ and let $\S(x): H\otimes H\rightarrow H\otimes
H\otimes \C[[x]]$ be a linear map with $\S(0)=1$. 
Assume that $V(H,\S)$ is of PBW type in the sense that
the $S(H\otimes t^{-1}\C[t^{-1}])$-module 
$\gr V(H,\S)=\coprod_{k\ge 0}V(H,\S)[k]/V(H,\S)[k-1]$,
obtained in Lemma \ref{lpbw-span-property},
is a free module.
Then there exists a unique weak quantum vertex structure on $V(H,\S)$ 
with ${\bf 1}$ as the vacuum vector such that
\begin{eqnarray}
Y(a(-1){\bf 1},x)=a(x)\;\;\;\mbox{ for }a\in H.
\end{eqnarray}
Furthermore, $V(H,\S)$ is a nondegenerate quantum vertex algebra.
\et

Notice that this theorem implies that 
if $V(H,\S)$ is of PBW type, $\S(x)$ is 
necessarily a unitary rational quantum Yang-Baxter operator.
We believe that if $\S(x)$ is a unitary rational
quantum Yang-Baxter operator,  $V(H,\S)$ is of PBW type.
Next we consider a special family, confirming our belief.

For the rest of this section we fix a finite-dimensional vector space $U$ and an
(invertible) element 
\begin{eqnarray}
R(x)=\sum_{n\ge 0}R_{n}x^{n}\in (\End U)[[x]]
\end{eqnarray}
such that $R_{i}R_{j}=R_{j}R_{i}$ for $i,j\ge 0$ with $R_{0}\;(=R(0))$ invertible.
Then we have $R(x_{1})R(x_{2})=R(x_{2})R(x_{1})$.

Define $R^{*}(x)\in (\End U^{*})[[x]]$ by
\begin{eqnarray}
\<R^{*}(x)u^{*},u\>=\<u^{*},R^{-1}(x)u\>
\end{eqnarray}
for $u\in U,\; u^{*}\in U^{*}$, where $\<\cdot,\cdot\>$ is the standard pairing between $U^{*}$ and $U$. Then
\begin{eqnarray}
\<R^{*}(x)u^{*},R(x)u\>=\<u^{*},u\>
\end{eqnarray}
for $u\in U,\; u^{*}\in U^{*}$.

Set
\begin{eqnarray}
H=U\oplus U^{*}
\end{eqnarray}
and equip $H$ with the skew-symmetric
bilinear form $\<\cdot,\cdot\>$ defined by
\begin{eqnarray}
\<u+u^{*},v+v^{*}\>=\<u^{*},v\>-\<v^{*},v\>
\end{eqnarray}
for $u,v\in U,\; u^{*},v^{*}\in U^{*}$.
Define $\S(x): H\otimes H\rightarrow H\otimes H\otimes \C[[x]]$ by
\begin{eqnarray}
& &\S(x)(u\otimes v)=R(-x)u\otimes R^{-1}(x)v,\\
& &\S(x)(u^{*}\otimes v^{*})=(R^{*})^{-1}(-x)u^{*}\otimes R^{*}(x)v^{*},\\
& &\S(x)(v^{*}\otimes u)=R^{*}(-x)v^{*}\otimes R(x)u,\\
& &\S(x)(u\otimes v^{*})=R^{-1}(-x)u\otimes (R^{*})^{-1}(x)v^{*}
\end{eqnarray}
for $u,v\in U,\; u^{*},v^{*}\in U^{*}$. 
We are going to show that if $R_{0}=1$, 
the universal vacuum $(H,\S)$-module $V(H,\S)$ 
is of PBW type. Our strategy is to use certain vertex operator algebras 
associated to Heisenberg Lie algebras.

\br{rbcsystem}
{\em Consider the Heisenberg Lie algebra
$\H=(U\oplus U^{*})\otimes \C[t,t^{-1}]\oplus \C c$,
where $c$ is central and 
\begin{eqnarray*}
& &[u\otimes t^{m},v\otimes t^{n}]=0
=[u^{*}\otimes t^{m},v^{*}\otimes t^{n}],\\
& &[u^{*}\otimes t^{m},u\otimes t^{n}]=\<u^{*},u\>\delta_{m+n+1,0}c
\end{eqnarray*}
for $u,v\in U,\; u^{*},v^{*}\in U^{*},\; m,n\in \Z$. 
In terms of generating function $a(x)=\sum_{n\in \Z}(a\otimes t^{n})x^{-n-1}$ for $a\in H$,
we have
\begin{eqnarray}
& &u(x_{1})v(x_{2})-v(x_{2})u(x_{1})=0,\\
& &u^{*}(x_{1})v^{*}(x_{2})-v^{*}(x_{2})u^{*}(x_{1})=0,\\
& &u^{*}(x_{1})v(x_{2})-v(x_{2})u^{*}(x_{1})
=\<u^{*},v\>x_{2}^{-1}\delta\left(\frac{x_{1}}{x_{2}}\right)
\end{eqnarray}
for $u,v\in U,\; u^{*},v^{*}\in U^{*}$. 
Set 
$$\H^{\ge 0}=(U\oplus U^{*})\otimes \C[t]\oplus \C c,\;\;\; 
\H^{-}=(U\oplus U^{*})\otimes t^{-1}\C[t^{-1}].$$ 
Consider $\C$ as an $\H^{\ge 0}$-module with $c$ acting as identity and with
$(U\oplus U^{*})\otimes \C[t]$ acting trivially. Form the induced $\H$-module
\begin{eqnarray}
V_{\H}=U(\H)\otimes _{U(\H^{\ge 0})}\C.
\end{eqnarray}
It is well known that $V_{\H}$ is an irreducible $\H$-module
and that $V_{\H}$ is naturally a simple conformal vertex algebra,
which we shall use.}
\er

The vertex algebra $V_{\H}$ satisfies the following universal property: 

\bl{luniversal}
Let $K$ be any vertex algebra and let $\psi: U\oplus U^{*}\rightarrow K$ 
be any linear map such that
\begin{eqnarray*}
& &[Y(\psi(u),x_{1}),Y(\psi(v),x_{2})]=0=[Y(\psi(u^{*}),x_{1}),Y(\psi(v^{*}),x_{2})],\\
& &[Y(\psi(u^{*}),x_{1}),Y(\psi(u),x_{2})]=\<u^{*},u\>x_{2}^{-1}\delta\left(\frac{x_{1}}{x_{2}}\right)
\end{eqnarray*}
for $u,v\in U,\; u^{*},v^{*}\in U^{*}$.
Then $\psi$ can be extended uniquely to a vertex algebra homomorphism from
$V_{\H}$ to $K$. 
\el

\begin{proof} From the assumption on $\psi$, $K$ is an $\H$-module
with $c$ acting as identity and with $a\otimes t^{n}$ acting as $\psi(a)_{n}$ 
for $a\in U\oplus U^{*},\; n\in \Z$. We have $c\cdot {\bf 1}={\bf 1}$ and
$\H^{\ge 0}{\bf 1}=0$.
Then there exists an $\H$-module homomorphism $f$ from $V_{\H}$ to $K$, 
sending the vacuum vector of $V_{\H}$ to the vacuum vector of $K$.
Since $U\oplus U^{*}$ generates $V_{\H}$ as a vertex algebra, it follows that 
$f$ is a vertex algebra homomorphism.
\end{proof}

Recall from Remark \ref{rLaurent} that Now, we have a vertex algebra 
$(\C[[t]],-d/dt)$ where
\begin{eqnarray}
Y(f(t),x)g(t)=\left(e^{-x(d/dt)}f(t)\right)g(t)=f(t-x)g(t)
\end{eqnarray}
for $f(t),g(t)\in \C[[t]]$.

\bl{laffine-endom}
There exist unique vertex algebra homomorphisms
$\Phi^{\pm}_{R}$ from $V_{\H}$ to the tensor product vertex algebra 
$\hat{V}_{\H}=V_{\H}\otimes \C[[t]]$ such that
\begin{eqnarray}
\Phi^{\pm}_{R}(u)=R^{\pm 1}(t)u,\ \ \ \ \
\Phi^{\pm}_{R}(u^{*})=(R^{*})^{\pm 1}(t)u^{*}
\end{eqnarray}
for $u\in U,\; u^{*}\in U^{*}$.
\el

\begin{proof} Since $U,U^{*}$ generate $V_{\H}$ as a vertex algebra,
the uniqueness is clear. Because of the obvious symmetry we
we shall only prove the existence of $\Phi^{+}_{R}$.
For $u,v\in U,\; u^{*},v^{*}\in U^{*}$, we have
\begin{eqnarray*}
& &[\hat{Y}(R(t)u,x_{1}),\hat{Y}(R(t)v,x_{2})]
=[Y(R(t-x_{1})u,x_{1}),Y(R(t-x_{2})v,x_{2})]=0,\\
& &[\hat{Y}(R^{*}(t)u^{*},x_{1}),\hat{Y}(R^{*}(t)v^{*},x_{2})]
=[Y(R^{*}(t-x_{1})u^{*},x_{1}),Y(R^{*}(t-x_{2})v^{*},x_{2})]=0,
\end{eqnarray*}
and 
\begin{eqnarray*}
& &[\hat{Y}(R^{*}(t)u^{*},x_{1}),\hat{Y}(R(t)u,x_{2})]\\
&=&[Y(R^{*}(t-x_{1})u^{*},x_{1}),Y(R(t-x_{2})u,x_{2})]\\
&=&\<R^{*}(t-x_{1})u^{*},R(t-x_{2})u\>
x_{2}^{-1}\delta\left(\frac{x_{1}}{x_{2}}\right)\\
&=&\<R^{*}(t-x_{1})u^{*},R(t-x_{1})u\>
x_{2}^{-1}\delta\left(\frac{x_{1}}{x_{2}}\right)\\
&=&\<u^{*},u\>x_{2}^{-1}\delta\left(\frac{x_{1}}{x_{2}}\right).
\end{eqnarray*}
In view of Lemma \ref{luniversal}, there exists a vertex algebra homomorphism 
$\Phi^{+}_{R}$ from $V_{\H}$ to $\hat{V}_{\H}$ with the desired properties.
\end{proof}

\bl{lpseudo-end}
Let $V$ be any vertex algebra. Equip $V\otimes \C[[t]]$ with the tensor product vertex algebra
structure, where $\C[[t]]$ is viewed as a vertex algebra as above.
A linear map $\Delta(t): V\rightarrow V\otimes \C[[t]]$ is a 
vertex algebra homomorphism if and only if the following relations hold
for $v\in V$:
\begin{eqnarray}
& &\Delta(t){\bf 1}={\bf 1}\;(={\bf 1}\otimes 1),\\
& &\Delta(t)Y(v,x)=Y(\Delta(t-x)v,x)\Delta(t).
\end{eqnarray}
\el

\begin{proof} Denote by $\hat{Y}$ the vertex operator map of the tensor product vertex algebra
$V\otimes \C[[t]]$. Note that for $a\in V, \;f(t)\in \C[[t]]$, we have
$$\hat{Y}(a\otimes f(t),x)=Y(a,x)f(t-x).$$
For $u,v\in V$, we have
\begin{eqnarray}
& &\hat{Y}(\Delta(t)u,x)\Delta(t)v=Y(\Delta(t-x)u,x)\Delta(t)v,\\
& &\Delta(t)\hat{Y}(u,x)v=\Delta(t)Y(u,x)v.
\end{eqnarray}
Then the assertion follows immediately.
\end{proof}

\bp{ppseudo-ext}
Assume that $R_{0}=R(0)=1$. There exist unique elements
$$\Delta_{R}^{\pm }(x)\in \Hom_{\C}(V_{\H},V_{\H}\otimes
\C[[x]])\subset (\End_{\C}V_{\H})[[x]]$$
such that
\begin{eqnarray}
& &\Delta_{R}^{\pm }(x)u=R^{\pm 1}(x)u,\ \ \ 
\Delta_{R}^{\pm}(x)u^{*}=(R^{*})^{\pm 1}(x)u^{*}
\;\;\;\mbox{ for }u\in U,\; u^{*}\in U^{*},\\
& &\Delta_{R}^{\pm }(x){\bf 1}={\bf 1},\\
& &\Delta^{\pm }_{R}(x)Y(v,x_{1})=Y(\Delta_{R}^{\pm}(x-x_{1})v,x_{1})\Delta^{\pm}_{R}(x)
\;\;\;\mbox{ for }v\in V_{\H}, \label{epseudo-eqn}\\
& &\Delta_{R}^{\pm}(x_{1})\Delta_{R}^{\pm
}(x_{2})=\Delta_{R}^{\pm}(x_{2})\Delta^{\pm }_{R}(x_{1}),
\ \ \ \ \Delta_{R}^{\pm}(x_{1})\Delta_{R}^{\mp
}(x_{2})=\Delta_{R}^{\mp}(x_{2})\Delta^{\pm }_{R}(x_{1}),\ \ \ \ 
\label{epseudo-eqn-comm}\\
& &\Delta_{R}^{\pm}(x)\Delta_{R}^{\mp}(x)=1.\label{einverse-relation}
\end{eqnarray}
For $u\in U,\; u^{*}\in U^{*}$, set
\begin{eqnarray}
& &u_{R}(x)=Y(u,x)\Delta^{+}_{R}(x)=u(x)\Delta^{+}_{R}(x),\\
& &u^{*}_{R}(x)=Y(u^{*},x)\Delta^{-}_{R}(x)=u^{*}(x)\Delta^{-}_{R}(x).
\end{eqnarray}
Then $V_{\H}$ is a vacuum $(H,\S)$-module of PBW type.
\ep

\begin{proof} In view of Lemma \ref{lpseudo-end}, such linear maps
$\Delta^{\pm}_{R}(t)$ if exist are vertex algebra homomorphisms (from $V_{\H}$
to $V_{\H}\otimes \C[[t]]$).
As $U\oplus U^{*}$ generates $V_{\H}$ as a vertex algebra,
the uniqueness is clear.
On the other hand, by Lemma \ref{laffine-endom}, 
we have vertex algebra homomorphisms
$\Phi^{\pm}_{R}$ from $V_{\H}$ to $V_{\H}\otimes \C[[t]]$ 
such that $\Phi^{\pm}_{R}(u)=R^{\pm 1}(t)u$ for $u\in U$ and
$\Phi^{\pm}_{R}(u^{*})=(R^{*})^{\pm 1}(t)u^{*}$ for $u^{*}\in U^{*}$. 
Let us alternatively denote the map $\Phi^{\pm}_{R}$ by $\Delta^{\pm}_{R}(t)$.
Then $\Delta^{\pm}_{R}(t){\bf 1}={\bf 1}$ and
the relation (\ref{epseudo-eqn}) holds.
For $u\in U$, we have
$$\Delta^{\pm}_{R}(x_{1})\Delta^{\pm}_{R}(x_{2})u
=R^{\pm 1}(x_{1})R^{\pm 1}(x_{2})u=R^{\pm 1}(x_{2})R^{\pm 1}(x_{1})u
=\Delta^{\pm}_{R}(x_{2})\Delta^{\pm}_{R}(x_{1})u.$$
{}From $R(x_{1})R(x_{2})=R(x_{2})R(x_{1})$ we have
$R^{*}(x_{1})R^{*}(x_{2})=R^{*}(x_{2})R^{*}(x_{1})$, so that
$$\Delta^{\pm}_{R}(x_{1})\Delta^{\pm}_{R}(x_{2})u^{*}
=(R^{*})^{\pm 1}(x_{1})R^{*}(x_{2})u^{*}=R^{*}(x_{2})R^{*}(x_{1})u^{*}
=\Delta^{\pm}_{R}(x_{2})\Delta^{\pm}_{R}(x_{1})u^{*}$$
for $u^{*}\in U^{*}$.
Suppose 
$$\Delta^{\pm}_{R}(x_{1})\Delta^{\pm }_{R}(x_{2})w=\Delta^{\pm}_{R}(x_{2})\Delta^{\pm}_{R}(x_{1})w$$
for some $w\in V_{\H}$.
For any $a\in U\oplus U^{*}$, using (\ref{epseudo-eqn}) we have
\begin{eqnarray*}
\Delta^{\pm}_{R}(x_{1})\Delta^{\pm}_{R}(x_{2})Y(a,x)w
&=&Y(\Delta^{\pm}_{R}(x_{1}-x)\Delta^{\pm}_{R}(x_{2}-x)a,x)\Delta^{\pm}_{R}(x_{1})\Delta^{\pm}_{R}(x_{2})w\\
&=&Y(\Delta^{\pm}_{R}(x_{2}-x)\Delta^{\pm}_{R}(x_{1}-x)a,x)\Delta^{\pm}_{R}(x_{2})\Delta^{\pm}_{R}(x_{1})w\\
&=&\Delta^{\pm}_{R}(x_{2})\Delta^{\pm}_{R}(x_{1})Y(a,x)w.
\end{eqnarray*}
Since $U\oplus U^{*}$ generates $V_{\H}$ as a vertex algebra, it follows that
$$\Delta^{\pm}_{R}(x_{1})\Delta^{\pm}_{R}(x_{2})w=\Delta^{\pm}_{R}(x_{2})\Delta^{\pm}_{R}(x_{1})w
\;\;\;\mbox{ for all }w\in V_{\H}. $$
This proves the first identity of (\ref{epseudo-eqn-comm}).
The second identity and (\ref{einverse-relation}) follow from the same argument. 
 
For $u,v\in U$, using (\ref{epseudo-eqn}) and (\ref{epseudo-eqn-comm}) we have
\begin{eqnarray*}
& &u_{R}(x_{1})v_{R}(x_{2})\\
&=&Y(u,x_{1})\Delta^{+}_{R}(x_{1})Y(v,x_{2})\Delta^{+}_{R}(x_{2})\\
&=&Y(u,x_{1})Y(\Delta^{+}_{R}(x_{1}-x_{2})v,x_{2})\Delta^{+}_{R}(x_{1})\Delta^{+}_{R}(x_{2})\\
&=&Y(\Delta^{+}_{R}(x_{1}-x_{2})v,x_{2})Y(u,x_{1})\Delta^{+}_{R}(x_{2})\Delta^{+}_{R}(x_{1})\\
&=&Y(R(x_{1}-x_{2})v,x_{2})\Delta^{+}_{R}(x_{2})Y(R^{-1}(x_{2}-x_{1})u,x_{1})
\Delta^{+}_{R}(x_{1}).
\end{eqnarray*}
Similarly, for $u^{*},v^{*}\in U^{*}$, we have
\begin{eqnarray*}
& &u^{*}_{R}(x_{1})v^{*}_{R}(x_{2})\\
&=&Y(u^{*},x_{1})\Delta_{R}^{-}(x_{1})Y(v^{*},x_{2})\Delta_{R}^{-}(x_{2})\\
&=&Y(u^{*},x_{1})Y(\Delta_{R}^{-}(x_{1}-x_{2})v^{*},x_{2})
\Delta_{R}^{-}(x_{1})\Delta_{R}^{-}(x_{2})\\
&=&Y(\Delta_{R}^{-}(x_{1}-x_{2})v^{*},x_{2})Y(u^{*},x_{1})
\Delta_{R}^{-}(x_{2})\Delta_{R}^{-}(x_{1})\\
&=&Y((R^{*})^{-1}(x_{1}-x_{2})v^{*},x_{2})\Delta_{R}^{-}(x_{2})
Y(R^{*}(x_{2}-x_{1})u^{*},x_{1})
\Delta_{R}^{-}(x_{1}).
\end{eqnarray*}
We also have for $u\in U,\;u^{*}\in U^{*}$
\begin{eqnarray*}
& &u_{R}(x_{1})v^{*}_{R}(x_{2})
-(R^{*}(x_{1}-x_{2})v^{*})_{R}(x_{2})(R(x_{2}-x_{1})u)_{R}(x_{1})\\
&=&Y(u,x_{1})\Delta_{R}(x_{1})Y(v^{*},x_{2})\Delta_{R}^{-}(x_{2})\\
& &-Y(\Delta^{+}_{R}(x_{1}-x_{2})v^{*},x_{2})\Delta_{R}^{-}(x_{2})
Y(\Delta^{+}_{R}(x_{2}-x_{1})u,x_{1})\Delta^{+}_{R}(x_{1})\\
&=&Y(u,x_{1})Y(\Delta^{+}_{R}(x_{1}-x_{2})v^{*},x_{2})\Delta^{+}_{R}(x_{1})\Delta_{R}^{-}(x_{2})\\
& &-Y(\Delta^{+}_{R}(x_{1}-x_{2})v^{*},x_{2})Y(u,x_{1})\Delta^{+}_{R}(x_{1})\Delta_{R}^{-}(x_{2})\\
&=&Y(u,x_{1})Y(R^{*}(x_{1}-x_{2})v^{*},x_{2})\Delta^{+}_{R}(x_{1})\Delta_{R}^{-}(x_{2})\\
& &-Y(R^{*}(x_{1}-x_{2})v^{*},x_{2})Y(u,x_{1})\Delta^{+}_{R}(x_{1})\Delta_{R}^{-}(x_{2})\\
&=&\<u,R^{*}(x_{1}-x_{2})v^{*}\>x_{2}^{-1}\delta\left(\frac{x_{1}}{x_{2}}\right)
\Delta^{+}_{R}(x_{1})\Delta_{R}^{-}(x_{2})\\
&=&\<u,R^{*}(0)v^{*}\>x_{2}^{-1}\delta\left(\frac{x_{1}}{x_{2}}\right)
\Delta^{+}_{R}(x_{2})\Delta_{R}^{-}(x_{2})\\
&=&\<u,v^{*}\>x_{2}^{-1}\delta\left(\frac{x_{1}}{x_{2}}\right).
\end{eqnarray*}
Substituting $v^{*}$ with $(R^{*})^{-1}(x_{1}-x_{2})v^{*}$ and
$u$ with $R^{-1}(x_{2}-x_{1})u$ and then exchanging $x_{1}$ with
$x_{2}$ we get
\begin{eqnarray}
& &(R^{-1}(x_{1}-x_{2})u)_{R}(x_{2})((R^{*})^{-1}(x_{2}-x_{1})v^{*})_{R}(x_{1})
-v^{*}_{R}(x_{1})u_{R}(x_{2})\nonumber\\
&=&\< R^{-1}(x_{1}-x_{2})u,(R^{*})^{-1}(x_{2}-x_{1})v^{*}\>
x_{1}^{-1}\delta\left(\frac{x_{2}}{x_{1}}\right)\nonumber\\
&=&-\<v^{*},R(x_{2}-x_{1}) R^{-1}(x_{1}-x_{2})u\>
x_{1}^{-1}\delta\left(\frac{x_{2}}{x_{1}}\right)\nonumber\\
&=&-\<v^{*},u\>x_{1}^{-1}\delta\left(\frac{x_{2}}{x_{1}}\right).
\end{eqnarray}
This proves that $V_{\H}$ is an $(H,\S)$-module. Clearly,
${\bf 1}$ is a vacuum-like vector.

On $V_{\H}$ we have an increasing filtration $\{E_{k}\}$ where for $k\ge 0$, 
$E_{k}$ is spanned by the vectors
$$ a^{(1)}(-m_{1})\cdots a^{(r)}(-m_{r}){\bf 1}$$
for $r\ge 0,\; a^{(i)}\in H=U\oplus U^{*},\; m_{i}\ge 1$ with
$m_{1}+\cdots +m_{r}\le k$. For $a\in H,\; m,n\in \Z$, we have
\begin{eqnarray*}
& &a(m)E_{n}\subset E_{n+m},\\
& &a(m)E_{n}\subset E_{n-m-1}\;\;\;\mbox{ for }m\ge 0.
\end{eqnarray*}
Furthermore, $\gr_{E}(V_{\H})$ is a 
free $S(H\otimes t^{-1}\C[t^{-1}])$-module.

For $a\in H$, set
$$a_{R}(x)=\sum_{n\in \Z}a_{R}(n)x^{-n-1}.$$
In the following we are going to prove that 
for each $k\ge 0$, $E_{k}$ is also linearly spanned by
the vectors
$$ a_{R}^{(1)}(-m_{1})\cdots a_{R}^{(r)}(-m_{r}){\bf 1}$$
for $r\ge 0,\; a^{(i)}\in H=U\oplus U^{*},\; m_{i}\ge 1$ with
$m_{1}+\cdots +m_{r}\le k$ and that
\begin{eqnarray}\label{ear-relation}
a_{R}(-m)w-a_{R}(-m)w\in E_{k+m-1}
\end{eqnarray}
for $a\in H,\; m\ge 1,\; w\in E_{k}$ with $k\ge 0$.
Then it will follow immediately  that $V_{\H}$ is a vacuum $(H,\S)$-module
cyclic on the vacuum-like vector ${\bf 1}$ and that
$V_{\H}$ is of PBW type.

Since $E_{0}=\C {\bf 1}$ and $\Delta^{\pm}_{R}(x){\bf 1}={\bf 1}$, we have 
$\Delta^{\pm}_{R}(x)w-w=0\in E_{-1}$ for $w\in E_{0}$.
For $a\in H=U\oplus U^{*}$ we have
$$\Delta^{\pm}_{R}(x)Y(a,x_{1})=Y(\Delta^{\pm}_{R}(x-x_{1})a,x_{1})\Delta^{\pm}_{R}(x).$$
Set $\Delta^{\pm}_{R}(x)a=\sum_{r\ge 0}a^{(r)}x^{r}$ with $a^{(0)}=a$.
For $m\in \Z$, we have
$$\Delta^{\pm}_{R}(x) a(m)
=a(m)\Delta^{\pm}_{R}(x)+\sum_{r\ge 1}\sum_{i\ge 0}\binom{r}{i}(-1)^{i} a^{(r)}(m+i)x^{r-i}\Delta^{\pm}_{R}(x).$$
It follows from induction on $k$ that for $w\in E_{k}$ with $k\ge 0$, we have 
\begin{eqnarray}
\Delta^{\pm}_{R}(x)w-w\in E_{k-1}[[x]]+xE_{k}[[x]].
\end{eqnarray}
As $a_{R}(x)=a(x)\Delta^{+}_{R}(x)$ for $a\in U$ and $a_{R}(x)=a(x)\Delta^{-}_{R}(x)$ 
for $a\in U^{*}$, we always have
\begin{eqnarray}
a_{R}(m)w-a(m)w\in \sum_{i\ge 0}a(m+i)E_{k-1}+\sum_{j\ge 1} a(m+j)E_{k}\subset E_{k-m-1},
\end{eqnarray}
proving (\ref{ear-relation}).
{}From this we have
\begin{eqnarray}
a_{R}^{(1)}(-m_{1})\cdots a_{R}^{(r)}(-m_{r}){\bf 1}-a^{(1)}(-m_{1})\cdots a^{(r)}(-m_{r}){\bf 1}
\in E_{m_{1}+\cdots +m_{r}-1}
\end{eqnarray}
for $r\ge 0,\; a^{(i)}\in H,\; m_{i}\ge 1$.
It follows from induction (on $k$) that $E_{k}$ 
is linearly spanned by the vectors
$$ a_{R}^{(1)}(-m_{1})\cdots a_{R}^{(r)}(-m_{r}){\bf 1}$$
for $r\ge 0,\; a^{(i)}\in H=U\oplus U^{*},\; m_{i}\ge 1$ with
$m_{1}+\cdots +m_{r}\le k$. Therefore $V_{\H}$ as an $(H,\S)$-module is of PBW type.
\end{proof}

Now we have:

\bc{cfinal}
Assume that $R_{0}=R(0)=1$.
The universal vacuum $(H,\S)$-module $V(H,\S)$ is of PBW type and
$V(H,\S)$ and $V_{\H}$ are isomorphic $(H,\S)$-modules.
\ec

\begin{proof} There exists a $T(H\otimes \C[t,t^{-1}])$-module homomorphism $f$
{}from $V(H,\S)$ onto
$V_{\H}$, sending the vacuum-like vector ${\bf 1}$ to ${\bf 1}$.
Clearly, $f$ reduces to an $S(H\otimes t^{-1}\C[t^{-1}])$-module homomorphism
$\bar{f}$ from $\gr V(H,\S)$ onto $\gr V_{\H}$.
Since $\gr V_{\H}$ is a free $S(H\otimes t^{-1}\C[t^{-1}])$-module
and since both  $\gr V_{\H}$ and $\gr V(H,\S)$ are cyclic
$S(H\otimes t^{-1}\C[t^{-1}])$-modules,  
$\bar{f}$ is an isomorphism.
Thus $\gr V(H,\S)$ is free. That is, $V(H,\S)$ is of PBW type.
{}From that $\bar{f}$ is an isomorphism, it follows that 
$f$ is an isomorphism.
\end{proof}

\end{document}